\documentclass[12pt]{amsart}

\usepackage{amsmath,amssymb,amsthm,epsfig }
\usepackage[all]{xy}
\usepackage{stmaryrd}
\usepackage{tikz}
\usepackage{pgfplots}
\pgfplotsset{compat=1.18}
\usepgfplotslibrary{fillbetween}
\usepackage{enumerate}
\usepackage{color}
\usepackage{bbm}
\usepackage{subcaption}
\usepackage{appendix}
\usepackage[margin=30truemm]{geometry}
\usepackage[colorinlistoftodos]{todonotes}
\usepackage{version}
\usepackage{enumitem}
\usepackage{setspace}
\usepackage{hyperref}
\usepackage[
  backend=biber,
  style=alphabetic,
  sorting=ynt
]{biblatex}
\usetikzlibrary{shapes.geometric, shapes.misc}
\DeclareFontFamily{OT1}{rsfs}{}
\DeclareFontShape{OT1}{rsfs}{n}{it}{<-> rsfs10}{}
\DeclareMathAlphabet{\curly}{OT1}{rsfs}{n}{it}

\newtheorem{Thm}{Theorem}[section]
\newtheorem{maintheorem}{Theorem}
\newtheorem{lem}[Thm]{Lemma}

\newtheorem{prop}[Thm]{Proposition}

\theoremstyle{remark}
\newtheorem{Rem}[Thm]{Remark}
\newtheorem{ex}[Thm]{Example}
\theoremstyle{definition}
\newtheorem{defn}[Thm]{Definition}
\newtheorem{propdef}[Thm]{Proposition-Definition}

\newtheoremstyle{algplain}
  {\topsep}               
  {\topsep}               
  {\itshape}              
  {}                      
  {\bfseries}             
  {.}                     
  {.5em}                  
  {\thmname{#1}\thmnumber{ #2\textsuperscript{$\diamond$}}
   \thmnote{ (#3)}}

\theoremstyle{algplain}
\newtheorem{algtheorem}[Thm]{Theorem}
\newtheorem{alglemma}[Thm]{Lemma}

\newtheorem{algcorollary}[Thm]{Corollary}

\theoremstyle{plain}

\newcommand{\Spec}{\mathop{\mathrm{Spec}}\nolimits}

\newcommand{\Hom}{\mathop{\mathrm{Hom}}\nolimits}
\newcommand{\Ext}{\mathop{\mathrm{Ext}}\nolimits}

\newcommand{\SL}{\mathop{\mathrm{SL}}\nolimits}
\newcommand{\GL}{\mathop{\mathrm{GL}}\nolimits}

\newcommand{\SU}{\mathop{\mathrm{SU}}\nolimits}

\newcommand{\st}{\mathop{\mathrm{st}}\nolimits}
\newcommand{\ord}{\mathop{\mathrm{ord}}\nolimits}
\newcommand{\gr}{\mathop{\mathrm{gr}}\nolimits}

\newcommand{\length}{\mathop{\mathrm{length}}\nolimits}

\newcommand{\wt}{\mathop{\mathrm{wt}}\nolimits}

\newcommand{\Sas}{{\mathop{\mathrm{Sas}}}}

\def\Lie{\operatorname{Lie}}
\def\GL{\operatorname{GL}}

\def\SL{\operatorname{SL}}

\def\End{\operatorname{End}}
\def\tr{\operatorname{tr}}

\def\rank{\operatorname{rank}}

\def\dim{\operatorname{dim}}

\def\det{\operatorname{det}}
\def\vol{\operatorname{vol}}
\def\deg{\operatorname{deg}}

\def\R{\mathbb{R}}

\hypersetup{hidelinks}

\excludeversion{memo}
\allowdisplaybreaks
\begin{document}

\title[Slope stability on Sasakian manifolds]{On the Slope Stability of vector bundles on Sasakian manifolds}
\author{Yohei Hada}

\maketitle
\thispagestyle{empty}
\maketitle
\begin{abstract}
    We study the variation problem of slope stability of $T$-linearized sheaves on a compact Sasakian manifold, or more generally, on an affine variety with a good $T$-action.
    In particular, we establish the wall-and-chamber structure for the slope stability of $T$-linearized torsion-free sheaves.
\end{abstract}
\tableofcontents
\section{Introduction}
For a compact Sasakian manifold $(X,g,\xi)$, the Reeb vector field can be varied while preserving the underlying CR structure, giving rise to the Reeb cone. 
It is therefore natural to ask how the slope stability condition for transverse holomorphic vector bundles, introduced in \cite{biswas2009vectorbundlessasakianmanifolds} and \cite{BARAGLIA2022108661}, changes under such variation. 
In this paper, we study this variation problem for a fixed $T$-linearized transverse holomorphic vector bundle $E$, where $T$ is the closure of the flow generated by a given Reeb vector field.
We vary the Reeb vector within $\mathfrak t={\rm Lie}(T)$, while keeping the underlying CR structure fixed. 
For each Reeb vector $\tau$ in the corresponding Reeb cone $\mathfrak t_+\subset \mathfrak t$, the $T$-linearization induces a canonical transverse holomorphic structure $E_\tau$, and we study how the stability of $E_\tau$ varies with $\tau$. 

The variation above has a natural algebraic description in terms of the associated affine cone to the given Sasakian manifold. The complex structures of the affine cone remain fixed as the Reeb vector varies in $\mathfrak t_+$. Pulling back $E$ to the K\"{a}hler cone and extending it across the vertex, we obtain a $T_\mathbb C$-linearized algebraic reflexive sheaf $\mathcal E$ on the affine cone $\mathcal C$. 
Hence the module of global sections
\[
M=\Gamma(\mathcal C,\mathcal E)
\]
is naturally graded by the character lattice $\mathbb X^*(T_\mathbb C)$, as is the coordinate ring $R=\Gamma(\mathcal C,\mathcal O_{\mathcal C})$. A Reeb vector field $\tau\in \mathfrak t_+$ induces an $\mathbb R$-grading by the homomorphism 
\[
\mathbb X^*(T_\mathbb C)\to \mathbb R\ ;\  \chi\mapsto \langle \chi,\tau\rangle.
\]
Thus the variation of the Reeb vector can be viewed algebraically as varying the $\mathbb R$-grading of a fixed $\mathbb X^*(T_\mathbb C)$-graded module $M$.

Thus, there are two natural stability notions corresponding to the Sasakian and the algebraic viewpoints respectively. On the Sasakian side, slope stability for transverse holomorphic vector bundles has been studied in the quasi-regular case by Biswas--Schumacher \cite{biswas2009vectorbundlessasakianmanifolds} and, more generally, by Baraglia--Hekmati \cite{BARAGLIA2022108661}.
The slope $\mu_{\Sas}$ in {\it op.cit} is defined in terms of the first basic Chern class and the transverse K\"{a}hler form, and stability is tested against saturated transverse holomorphic coherent subsheaves. On the algebraic side, we use the stability notion for finitely generated torsion-free $\mathbb R$-graded modules introduced in \cite{hada2026optimalalgebraictangentcone}. 
The slope $\mu_{\gr}$ in {\it op.cit} is defined in terms of the first two asymptotic Hilbert coefficients and stability is tested against saturated graded submodules. Although these two notions arise from rather different constructions, they are naturally related through the correspondence between transverse holomorphic bundles and graded modules on the affine cone.

Our first main result is that these two stability notions coincide.
\begin{maintheorem}[Theorem \ref{equality of slopes} and Theorem \ref{equivalence of stabilities}]
Let $(X,g,\xi)$ be a compact Sasakian manifold of dimension $2n+1$, with $n\geq 1$, and let $T$ be the compact torus defined as the closure of the Reeb flow generated by $\xi$.
For a $T$-linearized transverse holomorphic vector bundle $E$ on $X$ and the corresponding $\mathbb X^*(T_\mathbb C)$-graded module $M$, we have
\[
\mu_{\Sas}(E_\tau)=\frac{1}{2\pi}\mu_{\gr}(M_\tau),
\]
for any $\tau\in \mathfrak t_+$, where $M_\tau$ is the $\mathbb R$-graded module induced by $M$ and the Reeb vector $\tau\in \mathfrak t_+$.
Furthermore, $E_\tau$ is $\mu_{\Sas}$-(semi)stable if and only if $M_\tau$ is $\mu_{\gr}$-(semi)stable.
\end{maintheorem}

Having identified the two stability notions, we study how stability varies as the Reeb vector $\tau$ moves in $\mathfrak t_+$. Our first result in this direction is the openness of stability: if $E_\tau$ is stable, then $E_{\tau'}$ remains stable for all $\tau'\in \mathfrak t_+$ sufficiently close to $\tau$.
\begin{maintheorem}[Theorem \ref{Openness of Stability I}]
Let $(X,g,\xi)$ be a compact Sasakian manifold of dimension $2n+1$ with $n\geq 1$, and $E$ be a $T$-linearized transverse holomorphic vector bundle. Assume that $E_\tau$ is stable for some $\tau\in \mathfrak t_+$. Then there is an open neighbourhood $U\subset \mathfrak t_+$ of $\tau$ such that $E_{\tau'}$ is stable for any $\tau' \in U$.
\end{maintheorem}

\begin{Rem}
In the setting of \cite[Section 3.5]{HL24}, openness of
transverse slope stability can also be deduced from
\cite[Proposition 3.6]{HL24}, under the simplicity
hypothesis imposed there.
\end{Rem}
We first give a differential geometric proof of this openness result. The main ingredients are a uniform spectral gap for the slopes of proper saturated subsheaves, obtained using the transverse Kobayashi--Hitchin correspondence proved in \cite{BARAGLIA2022108661}, and the derivative estimate for slope function, obtained from the theory of index characters developed by \cite{CS18}.

While the differential geometric argument is local in the Reeb cone, the algebraic approach leads to a global finiteness result. More precisely, we show that only finitely many saturated submodules can occur as maximal destabilizing submodules as the Reeb vector varies. 
An important observation is that, for any $\tau\in \mathfrak t_+$, the maximal destabilizing submodule of $M_\tau$ is induced by an $\mathbb X^*(T_\mathbb C)$-graded submodule of the fixed module $M$. Equivalently, the maximal destabilizing subsheaf of $E_\tau$ is $T$-invariant. 
As a consequence, the stable and semistable loci are determined by finitely many slope inequalities, and the Reeb cone admits a finite Harder--Narasimhan chamber decomposition.
\begin{maintheorem}[Theorem \ref{finiteness of walls}]\label{Theorem C}
Let $\mathcal C$ be a normal affine variety with good $T_\mathbb C$-action, and $M$ be a $\mathbb X^*(T_\mathbb C)$-graded torsion-free $R=\Gamma(\mathcal C,\mathcal O_\mathcal C)$-module. 
Define the stable cone $\mathfrak t_+^{\st}(M)$ and semistable cone $\mathfrak t_+^{{\rm ss}}(M)$ as follows:
    \[
    \mathfrak t_+^{\st}(M):=\{\tau\in \mathfrak t_+\mid M_\tau\text{ is stable}\},\ \mathfrak t_+^{{\rm ss}}(M):=\{\tau\mid \text{$M_\tau$ is semistable}\}.
    \]
    Then, the following properties hold:
    \begin{enumerate}
        \item $\mathfrak t_+^{\st}(M)$ is an open cone.
        \item There are finitely many saturated graded submodules $N_1,\dots,N_p\subset M$ such that 
        \[
        \mathfrak t_+^{{\rm ss}}(M)=\bigcap_{i=1}^p \{\tau\in \mathfrak t_+\mid \mu(N_{i,\tau})\leq \mu(M_\tau)\}.
        \] 
        \item There is a canonical decomposition of $\mathfrak t_+$ into finitely many chambers $\sigma_1,\dots,\sigma_K$ such that 
        \begin{itemize}
            \item For each $\sigma_i$, there is a filtration $F^i_\bullet: 0\subsetneq M^i_1\subsetneq M^i_2\subsetneq \dots\subsetneq M^i_{l_i}=M$ by saturated $\mathbb X^*(T_\mathbb C)$-graded submodules such that $F^i_{\bullet}$ gives the Harder-Narasimhan filtration of $M_\tau$ for every $\tau \in \sigma_i$.
            \item If $i\neq j$, the corresponding Harder-Narasimhan filtrations $F^i_\bullet$ and $F^j_\bullet$ are different.
        \end{itemize}
    \end{enumerate}
\end{maintheorem}
The proof of Theorem \ref{Theorem C} is based on reducing the possible slope functions to a finite-dimensional combinatorial problem. For a saturated submodule $N\subset M$ of fixed rank $s$, its determinant module can be embedded into one of finitely many shifts of $R$, and its slope can be expressed in terms of codimension-one multiplicities of the cokernel. 
Using the Altmann--Hausen description of the affine $T_\mathbb C$-varieties \cite{AH06}, and the Petersen--S\"{u}{\ss} description of their invariant prime divisors \cite{PS11}, we show that the contributions of all $T_\mathbb C$-invariant prime divisors can be written as $\mathbb Z_{\geq 0}$-linear combinations of finitely many positive functions on the Reeb cone. 
Thus, each slope function is encoded by a vector in $\mathbb Z_{\geq 0}^L$ for some $L\in \mathbb Z_{>0}$, in such a way that the coordinatewise order reverses the order of slopes. 
Dickson's lemma then implies that only finitely many maximal slope functions can occur in each rank. This, together with the uniqueness of the maximal destabilizing submodules, gives the finiteness of possible maximal destabilizing submodules, and iterating the argument on successive quotients yields the finite Harder--Narasimhan chamber decomposition.

Variation of slope stability with respect to geometric parameters has been studied in several related settings. Delloque \cite{Delloque25} studied stability of torsion-free sheaves on balanced manifolds under variation of the balanced class, where the stable locus is an open convex cone locally determined by finitely many linear inequalities. 
Harder--Narasimhan chamber structures with respect to movable curve classes were studied by Neumann \cite{neumann2009decomposition}. 
In contrast to the usual variation of polarization on a fixed complex space, varying the Reeb vector changes the transverse space itself. For a quasi-regular Reeb vector $\tau$, this transverse space is the orbifold quotient $Y_\tau=X/U(1)_\tau$, which therefore varies with $\tau$.

We conclude with several explicit examples. We first study examples on $\mathbb A^3$ with the standard diagonal $(\mathbb C^*)^3$-action. Since the affine cone is factorial, all slope functions are necessarily linear on the Reeb cone. We then turn to the conifold, which is not $\mathbb Q$-factorial, and exhibit nonlinear walls.

The paper is organized as follows. In Section 2, we recall transverse holomorphic vector bundles on Sasakian manifolds and the corresponding notions of slope stability, and the algebraic stability theory for graded modules. 
In Section 3, we establish a correspondence between transverse slope stability and graded module stability. 
Section 4 is devoted to the variation problem: in Section 4.1 we prove the openness of stability by differential-geometric methods, while in Section 4.2 we prove that stability is governed by a wall-and-chamber structures with finitely many walls by algebro-geometric methods. 
In Section 5, we discuss explicit examples, including weighted projective stacks and the conifold.
\section*{Notation}
Throughout this paper, a variety means a separated integral scheme of finite type over $\mathbb C$. We use the following notation:
\begin{itemize}
    \item $X$: a compact connected smooth manifold of dimension $2n+1$, where $n\geq 1$.
    \item $T$: a compact torus.
    \item $\mathfrak t$: the Lie algebra of $T$.
    \item $T_\mathbb C$: the complexification of $T$; in the purely algebraic setting, $T_\mathbb C$ denotes an algebraic torus.
    \item $\mathbb X^*(T_\mathbb C)$: the character lattice of $T_\mathbb C$.
    \item $\mathcal C$: the affine cone associated with a Sasakian manifold $X$; in the purely algebraic setting, a normal affine variety with a good $T_\mathbb C$-action.
    \item $\mathcal C^*=\mathcal C\setminus \{0\}$: the punctured affine cone of a Sasakian manifold $X$.
    \item $R=\Gamma(\mathcal C,\mathcal O_\mathcal C)$: the coordinate ring of $\mathcal C$.
    \item $M$: a finitely generated torsion-free $\mathbb X^*(T_\mathbb C)$-graded $R$-module.
\end{itemize}
\section{Preliminary}
\subsection{Vector bundles on Sasakian manifolds}
Let $(X,g,\xi)$ be a smooth Sasakian manifold of dimension $2n+1$, i.e. satisfying the following properties:
\begin{enumerate}
    \item $\xi$ is Killing vector field with $\|\xi\|_g\equiv 1$, i.e. $\mathcal L_\xi g=0$.
    \item Let $\Phi\in A^0({\rm End}(T_X))$ be $\nabla\xi$, where $\nabla$ is the Levi-Civita connection on $(X,g)$. Then we have 
    \[
    (\nabla_v\Phi)(w)=g(\xi,w)v-g(v,w)\xi
    \]
    for any $x\in X$ and $v,w\in T_{X,x}$.
\end{enumerate}
It is well known that these conditions are equivalent to the cone $(X\times \mathbb R_{>0}, (dr)^2+r^2g)$ being K\"{a}hler, where complex structure $J$ is given by $J|_F=\Phi|_F$, $J(r\partial_r)=\xi$, and $J\xi=-r\partial_r$, where $F=\xi^\perp\subset T_X$. We call $\xi$ a Reeb vector of $(X,g)$. 

For a Sasakian manifold $(X,g,\xi)$, we have canonical map
\[
q: T_X\to T_X/F\simeq \mathcal F_\xi:=\langle \xi\rangle,
\]
This can be written as $q=\eta\otimes \xi$ for some $\eta\in A^1(X)$. This $\eta$ is a contact form of $X$, i.e. $(d\eta)^n\wedge \eta$ is a volume form on $X$. Indeed, we have the following lemma:
\begin{lem}\label{almost Kahler transverse structure}
    For any $v,w\in F_x$ and $x\in X$, we have
    \[
    d\eta(v,w)=2g(\Phi v,w)
    \]
\end{lem}
\begin{proof}
    Let $\tilde v, \tilde w$ be sections of $F$ with $\tilde v(x)=v$, $\tilde w(x)=w$. Then we have 
    \begin{align*}
        d\eta(\tilde v,\tilde w)=&\tilde v(\eta(\tilde w))-\tilde w(\eta(\tilde v))-\eta([\tilde v,\tilde w])\\
        =&-\eta([\tilde v,\tilde w])\\
        =&-g(\xi,\nabla_{\tilde v}\tilde w-\nabla_{\tilde w}\tilde v)\\
        =&g(\Phi(\tilde v),\tilde w)-g(\tilde v,\Phi(\tilde w))\\
        =&2g(\Phi\tilde v,\tilde w).
    \end{align*}
    Thus we get the lemma.
\end{proof}
Thus $\Phi|_F$ also defines a transverse almost Kähler structure on the Reeb foliation $\mathcal F_\xi$. Indeed, this transverse almost Kähler structure is integrable and thus $F_\xi$ is a Kähler foliation. 

Since ${\rm Isom}(X,g)$ is compact, the closure of the Reeb flow $\phi: \mathbb R\to {\rm Isom}(X,g)$ is a torus $T\simeq (S^1)^r$. There is a unique element of $\mathfrak t:=\Lie(T)$, which is also denoted by $\xi$, such that the induced vector field $\left.\frac{d}{dt}\right|_{t=0}\exp(t\xi)$ coincides with $\xi$. We call this $\xi$ the Reeb vector. Let $\Lambda:=\ker(\exp: \Lie(T)\to T)$. Then $\Lambda\simeq \mathbb Z^r$ and this induces a natural lattice structure on $\mathfrak t$. If $\lambda \in \mathfrak t$ is in $\mathfrak t_\mathbb Q:=\mathbb Q\cdot\Lambda$, then every orbit of the flow generated by $\lambda$ is closed.

We may view $(X,g,\xi)$ as an $\mathbb R$-equivariant manifold with transverse Kähler structure via the Reeb flow $\phi$. Thus, it is natural to consider $\mathbb R$-equivariant structures on $X$. Furthermore, each transverse slice has a natural Kähler structure induced by the transverse Kähler structure on $X$. For example, we can define the following notion:
\begin{defn}
    A complex differential form $\alpha$ on a Sasakian manifold $(X,g,\xi)$ is called basic if it satisfies the following equations:
    \[
    i_\xi\alpha =0, \ i_\xi d\alpha=0
    \]
    For any open set $U\subset X$, we denote by $\Omega^p_B(U)$ the space of all basic $p$-forms.
\end{defn}
Note that $\alpha$ being basic is equivalent to $i_\xi \alpha=L_\xi\alpha=0$, which means that $\alpha$ descends to a differential form on the quotient differential stack $[X/\mathbb R]$. Thus we have the following proposition:
\begin{prop}
    Let $(X,g,\xi)$ be a Sasakian manifold. Then we have the following Hodge-type decomposition:
    \[
    \Omega^p_B(X)=\bigoplus_{i=0}^p\Omega_B^{i,p-i}(X),
    \]
    and $d\eta\in \Omega^{1,1}_B(X)$.
\end{prop}
\begin{proof}
    We can decompose $F\otimes \mathbb C=F^{1,0}\oplus F^{0,1}$, where $F^{1,0}$, $F^{0,1}$ are the eigenspaces of $\Phi|_F$ with eigenvalues $\sqrt{-1},-\sqrt{-1}$, respectively. Let $w$ be a smooth section of $F\otimes \mathbb C$. Then we have 
    \begin{equation}\label{Reeb vector field preserves CR structure}
        (\nabla_\xi\Phi)(w)=(\mathcal L_\xi \Phi)(w)=0,
    \end{equation}
    Thus $\Phi$ preserves $F^{1,0}$ and $F^{0,1}$. Indeed, for $w\in \Gamma(F^{1,0})$, we have
    \[
    \Phi(w)=\sqrt{-1}w,
    \]
    Thus we have
    \[
    0=(\mathcal L_\xi\Phi)(w)=[\xi,\Phi(w)]-\Phi([\xi,w])=\sqrt{-1}[\xi,w]-\Phi([\xi,w]).
    \]
    Thus the Reeb flow preserves $F^{1,0}$, and similarly $F^{0,1}$. Thus it also preserves $\Omega^{1,0}_B$ and $\Omega^{0,1}_B$ by the equality
    \[
    \mathcal L_\xi(\alpha(w))=\alpha([\xi,w])+(\mathcal L_\xi\alpha)(w).
    \]
    Take $\alpha\in \Omega_B^1(X)$. If $\mathcal L_\xi\alpha=0$, then $0=\mathcal L_\xi\alpha^{1,0}=\mathcal L_\xi\alpha^{0,1}$ since the Reeb flow preserves the decomposition $\Omega^1_B(X)=\Omega^{1,0}_B(X)\oplus \Omega^{0,1}_B(X)$, i.e. $\alpha^{1,0}$ and $\alpha^{0,1}$ are basic since it is trivial that $i_\xi\alpha^{0,1}=i_\xi\alpha^{1,0}=0$. Thus we get the first statement of the proposition by taking exterior products. The last statement follows from Lemma \ref{almost Kahler transverse structure}, and $\mathcal L_\xi \eta=0$, $\eta(\xi)=1$.
\end{proof}
We consider the basic de-Rham cohomology $H^p_B(X,\mathbb C)$ obtained by the following complex: 
\[
0\to \Omega^0_B(X)\overset{d}\to \Omega^1_B(X)\overset d\to \cdots.
\]
Since the Reeb foliation of a compact Sasakian manifold is homologically orientable and transversely K\"{a}hler, there is a corresponding Hodge decomposition (see: \cite[Theorem 3.4.6]{ElKacimi1990})
\[
H^p_B(X,\mathbb C)=\bigoplus_{i+j=p}H^{i,j}_B(X).
\]
We next recall the definition of transverse vector bundles on Sasakian manifolds.
\begin{defn}\cite[Section 3.2]{biswas2009vectorbundlessasakianmanifolds}
    Let $E$ be a complex vector bundle on Sasakian manifold $(X,g,\xi)$. Let $S\subset T_X\otimes \mathbb C$ be a sub-bundle which is integrable, i.e. $[S,S]\subset S$. A partial connection on $E$ in the direction of $S$ is a first-order differential operator
    \[
    D_S: \Gamma(E)\to \Gamma(S^*\otimes E)
    \]
    which satisfies the following Leibniz rule:
    \[
    D_S(fs)=fD_S(s)+q_S(df)\otimes s,
    \]
    where $q_S: A^1(X)_\mathbb C\to \Gamma(S^\lor)$ is induced by the dual of the inclusion map $S\subset T_X\otimes \mathbb C$. $D_S$ extends to $\Gamma(\wedge^pS^*\otimes E)\to \Gamma(\wedge^{p+1}S^*\otimes E)$ in a natural way since $S$ is integrable. The curvature of partial connection $D_S$ is defined by $R(D_S)=D_S^2: \Gamma(E)\to \Gamma(\wedge^2S^*\otimes E)$, which is in fact a bundle homomorphism. We call $D_S$ flat if $R(D_S)=0$. 
\end{defn}
\begin{defn}\cite[Definition 3.2, Definition 3.3]{biswas2009vectorbundlessasakianmanifolds}
    A transverse complex vector bundle on Sasakian manifold $(X,g,\xi)$ is a pair $(E,D_\xi)$ consisting of complex vector bundle $E$ and a partial connection $D_\xi$ on $E$ in the direction of $F_\xi$. A holomorphic structure on $(E,D_\xi)$ is a flat partial connection $D''$ on $E$ in the direction of $F_\xi\oplus F^{0,1}$, which is in fact integrable, such that restriction of $D''$ to $F_\xi$ coincides with $D_\xi$. $s\in \Gamma(E)$ is called transverse holomorphic if $D''s=0$. 
    A transverse holomorphic local frame of $E$ is a local frame $(s_i)_{i=1}^r$ consisting of transverse holomorphic sections.
\end{defn}
Note that $D_\xi$ is automatically flat since $F_\xi$ has rank $1$, and hence it defines parallel transport along the Reeb trajectories. If we have $\mathbb R$-invariant hermitian metric on transverse holomorphic vector bundle $(E,D_\xi, D'')$, which means that hermitian metric on $E$ with $D_\xi h=0$, then we have associated connection $D=D'+D''$, where $D': \Gamma(E)\to \Gamma((F^{1,0})^*\otimes E)$ on $E$ by 
\[
d''(h(s_i,s_j))=h(s_i,D's_j)
\]
for any transverse holomorphic local frame $(s_i)_{i=1}^r$ on $E$, where $d=d'+d''+d_\xi$ is the decomposition of $d$ with respect to $T_X=F^{1,0}\oplus F^{0,1}\oplus F_\xi$. Note also that $d_\xi(h(s_i,s_j))=0$ since $D_\xi h=0$ and $D_\xi s_i=0$.

Since $D$ is compatible with the flat partial connection $D_\xi$ and is invariant under the Reeb flow, its curvature satisfies 
\[
i_\xi R(D)=0, \ \mathcal L_\xi R(D)=0.
\]
Moreover, $R(D)$ is transverse of type $(1,1)$. Thus, we can define the Chern forms by the following:
\[
\det\left(1+\frac{\sqrt{-1}}{2\pi}R(D)\right)=1+\sum_{i=1}^n c_i(E,h).
\]
By the usual Chern-Weil argument, $c_i(E,h)$ are closed basic forms and the corresponding cohomology class $[c_i(E,h)]\in H^{i,i}_B(X)$ does not depend on the choice of $h$. This defines the transverse Chern class $c_i(E)\in H^{i,i}_B(X)$. 

Transverse coherent sheaves on a Sasakian manifold $(X,g,\xi)$ is defined as follows. 
\begin{defn}\cite[Definition 3.7]{biswas2009vectorbundlessasakianmanifolds}
    Define the sheaf $\mathcal O_X$ of transverse holomorphic functions to be 
    \[
    \mathcal O_X(U)=\{f\in \Omega_B^0(U)\mid d''f=0\}.
    \]
    It is straightforward to show that $\mathcal O_X$ is a sheaf of rings. 
    A transverse coherent sheaf on $X$ is an $\mathcal O_X$-module $\mathcal E$ which is locally finitely presented; namely, for every $x\in X$, there is an open neighbourhood $U_x$ of $x$, integers $a_x,b_x\geq 0$, and a presentation
    \[
    \mathcal O_{U_x}^{b_x}\to \mathcal O_{U_x}^{a_x}\to \mathcal E|_{U_x}.
    \]
\end{defn}
Let $E$ be a transverse holomorphic vector bundle, then we can make a corresponding locally free $\mathcal E$ of $\mathcal O_X$-module, conversely, we can make a transverse holomorphic vector bundle from locally free sheaf of $\mathcal O_X$-module.

Let $\mathcal E$ be a coherent $\mathcal O_X$-module. For each point $x\in X$, we can take a foliated chart $\phi: X\supset U\simeq V\times (-\epsilon,\epsilon)$, where $U$ is an open neighbourhood of $x$, $0\in V\subset \mathbb C^n$ is an open set, and $\phi(x)=(0,0)$, which also reflects its transverse holomorphic structure. Fix such $\phi$ and we identify $U$ and $V\times (-\epsilon,\epsilon)$. Since transverse holomorphic functions are invariant under the Reeb flow, we have 
\[
\mathcal E|_U=\pi_1^*\mathcal E|_{V\times \{0\}},
\]
where $\pi_1: V\times (-\epsilon,\epsilon)\to V$ is the first projection. Thus, locally, the coherent sheaf is the same as the coherent sheaf on $V$ and thus we can define local properties and local constructions of coherent sheaves on $X$ by taking local chart. In particular, we can define torsion-freeness, reflexivity, and determinant line bundle of the transverse holomorphic coherent sheaves.
The stability condition in this setting was defined in \cite{biswas2009vectorbundlessasakianmanifolds} and \cite{BARAGLIA2022108661}.

In the quasi-regular setting considered by Biswas and Schumacher \cite{biswas2009vectorbundlessasakianmanifolds}, transverse holomorphic bundles and coherent sheaves descend to the orbifold quotient, so their determinants admit transverse Hermitian metrics. For a general irregular Reeb foliation, however, the existence of such metrics is not automatic. This issue was addressed by Baraglia–Hekmati \cite{BARAGLIA2022108661} using foliated resolutions of singularities. In the present $T$-linearized setting, as we see below, every saturated transverse holomorphic subsheaf of a $T$-linearized transverse holomorphic vector bundle $E$ is again $T$-linearized. Therefore, averaging an arbitrary hermitian metric on determinant line bundle over the compact torus provides a required transverse Hermitian metric on determinant line bundle, allowing us to use the simpler determinant definition of the degree.

\begin{defn}
    A $T$-linearized transverse holomorphic vector bundle $E$ on Sasakian manifold $(X,g,\xi)$ is transverse holomorphic vector bundle $E$ such that there is a lift of transverse holomorphic $T$-action to $E$ such that 
    \[
    D_\xi s=\mathcal L_\xi s\otimes \eta,
    \]
    where 
    \[
    \mathcal L_\xi s(x)=\lim_{t\to 0}\frac{\tilde \phi_t^*s(\phi_t(x))-s(x)}{t}\in E_x,\ \tilde \phi_t^*: \phi_t^*E\simeq E; s\in E_{\phi_t(x)}\mapsto \exp(-\xi t)\cdot s\in E_x.
    \]
    For a transverse torsion-free coherent sheaf $\mathcal E$ on $(X,g,\xi)$, we say that $\mathcal E$ is $T$-linearized if there is a lift of transverse holomorphic $T$-action to $\mathcal E$ such that the induced action on its locally free locus satisfies the compatibility condition with the Reeb flow as above. 
\end{defn}

\begin{lem}\label{saturated subsheaves of $T$-linearized bundle are $T$-invariant}
    Let $(X,g,\xi)$ be a compact Sasakian manifold, and $\mathcal E$ be a $T$-linearized transverse torsion-free coherent sheaf on $(X,g,\xi)$. Let $\mathcal F$ be a saturated transverse coherent subsheaf of $\mathcal E$. Then $\mathcal F$ is a $T$-invariant subsheaf. Moreover, $\det(\mathcal F)$ has a natural $T$-linearization.
\end{lem}
\begin{proof}
    Let $U$ be the locally free locus of $\mathcal E$, and $j: U\hookrightarrow X$ be the inclusion map. Since $\mathcal F=\mathcal E\cap j_*(\mathcal F|_U)$, it suffices to prove that $\mathcal F|_U$ is $T$-invariant. 
    First, note that $\exp(t\xi)\mathcal F\subset \mathcal F$ since transverse coherent subsheaves of $\mathcal E$ are preserved by the Reeb flow. Let $s$ be a local section of $\mathcal F$ on an open set $V\subset U$. For any $g\in T$, there exists a sequence $t_j$ such that 
    \[
    \lim_{j\to \infty}\exp(t_j\xi)=g.
    \]
    Take an arbitrary $x\in V$ and a small open neighbourhood $V'$ of $gx$ with $\bar V'\subset gV$. Then $V'\subset \exp(t_j\xi)V$ for $j\gg 0$. Since the action of $T$ on $\mathcal E|_U$ is smooth, we have 
    \[
    \lim_{j\to \infty} \exp(t_j\xi)s|_{V'}=gs|_{V'},
    \]
    in the $C^\infty$-Fréchet topology of $E(V)$. Since $\mathcal F$ is saturated, there is a morphism of transverse holomorphic vector bundles $\phi_{V'}: E|_{V'}\to E'_{V'}$ on $V'$ such that $\ker(\phi_{V'})=F|_{V'}$ after shrinking $V'$ if necessary. 
    Indeed, after shrinking $V'$ if necessary, the torsion-free sheaf $\mathcal Q_{V'}=(E/\mathcal F)|_{V'}$ admits an embedding $\mathcal Q_{V'}\to E'_{V'}$ into a locally free sheaf $E'_{V'}$.
    This shows that $\mathcal F(V')$ is a closed subspace of $E(V')$ with respect to $C^\infty$-Fréchet topology. Thus $gs|_{V'}\in F(V')$, and thus $gs\in \mathcal F(gV)$ since $x\in V$ is arbitrary. 

    Since the operations of taking exterior powers and reflexive hulls are functorial, the $T$-linearization of $\mathcal F$ induces a $T$-linearization of $\det(\mathcal F)=(\wedge^{\rank \mathcal F}\mathcal F)^{**}$.
\end{proof}

\begin{defn}\label{BS-stability}
    Let $(X,g,\xi)$ be a compact Sasakian manifold of dimension $2n+1$. For a $T$-linearized transverse torsion-free coherent sheaf $\mathcal E$ on $X$, define the degree of $\mathcal E$ as 
    \[
    \deg_{\Sas} (\mathcal E)=\frac{1}{\vol(X)}\int_Xc_1(\det \mathcal E)\wedge (d\eta)^{n-1}\wedge \eta,
    \]
    where 
    \[
    \vol(X)=\int_X(d\eta)^n\wedge \eta.
    \]
    The slope $\mu_{\Sas}(\mathcal E)$ of $\mathcal E$ is defined as 
    \[
    \mu_{\Sas}(\mathcal E)=\frac{\deg_{\Sas}(\mathcal E)}{\rank \mathcal E}
    \]
    The $T$-linearized transverse torsion-free coherent sheaf $\mathcal E$ is semistable if $\mu_{\Sas}(\mathcal F)\leq \mu_{\Sas}(\mathcal E)$ for any saturated transverse coherent submodule $\mathcal F\subset \mathcal  E$ with $0<\rank \mathcal F<\rank \mathcal E$, stable if the inequality is strict for any $\mathcal F$ as above, and polystable if $\mathcal E$ is direct sum of stable sheaves with the same slope.
\end{defn}

\begin{Rem}
    By \cite[Proposition 3.21]{BARAGLIA2022108661}, this definition coincides with the general definition given by Baraglia--Hekmati.
\end{Rem}
There is a Kobayashi-Hitchin correspondence in this setting, which is indicated by Biswas--Schumacher \cite{biswas2009vectorbundlessasakianmanifolds} in the quasi-regular case and proved by Baraglia--Hekmati \cite{BARAGLIA2022108661} in more general settings:
\begin{Thm}
    Let $(X,g,\xi)$ be a compact Sasakian manifold and $E$ be a $T$-linearized transverse holomorphic vector bundle on $X$. Then the following conditions are equivalent:
    \begin{enumerate}
        \item $E$ is polystable.
        \item $E$ admits a $T$-invariant transverse Hermitian--Einstein metric $h$; namely, $h$ is a $T$-invariant transverse Hermitian metric such that there is a constant $\lambda\in \mathbb R$ with 
        \[
        \sqrt{-1}\Lambda R(E,h)=\lambda \cdot I_E,
        \]
        where $\Lambda$ is the formal adjoint of $L: \Omega^0_B\to \Omega^{1,1}_B; f\mapsto f\otimes \left(\frac{1}{2}d\eta\right)$ with respect to $g$.
    \end{enumerate}
\end{Thm}
\begin{proof}
    By \cite{BARAGLIA2022108661}, we know that the existence of transverse Hermitian--Einstein metric is equivalent to the polystability of $E$. 
    Thus, it suffices to show that $E$ admits $T$-invariant transverse Hermitian--Einstein metric if $E$ is polystable. In this case, we know that there is transverse Hermitian--Einstein metric $h$ on $E$. 
    Let $g\in T$ be an arbitrary element. Then the Chern connection of the metric $g^*h$ is $g^{-1}\circ D\circ g$, where $D$ is the Chern connection of $h$. Thus we have $R(g^*h)=g^{-1}\circ R(h)\circ g$. 
    This and the assumption that $h$ is transverse Hermitian--Einstein show that $g^*h$ is also transverse Hermitian--Einstein. However, by the uniqueness of the Chern connection of transverse Hermitian--Einstein metric, which is proved by the same argument as the proof of \cite[Proposition 4.5]{BARAGLIA2022108661}, the Chern connection associated to two transverse Hermitian--Einstein metrics on the same transverse holomorphic vector bundle coincide. Thus, we have $g^{-1}\circ D\circ g=D$. In particular, $g^*h$ is $D$-parallel.

    Let $\mu$ be the normalized Haar measure of $T$ and let $h_T$ be the metric defined by 
    \[
    h_T=\int_T g^*hd\mu(g).
    \]
    Then $h_T$ is also $D$-parallel and thus the Chern connection of $h_T$ is also $D$. Thus $h_T$ is a $T$-invariant transverse Hermitian--Einstein metric. 
\end{proof}
\subsection{Algebraic stability}
From now on, $(X,g,\xi)$ is a compact Sasakian manifold of dimension $2n+1$ with $n\geq 1$. Define the Reeb cone $\mathfrak t_+$ by 
\begin{equation}\label{Reeb cone}
\mathfrak t_+:=\{\tau\in \mathfrak t\mid \eta(\tau)>0\},
\end{equation}
where $\eta(\tau)> 0$ means that $\eta(\tau_x)>0$ for any $x\in X$. Then, for every $\tau \in \mathfrak t_+$, we obtain a Sasakian structure on $X$ with the same CR structure as the original one and with Reeb vector $\tau$. Indeed, we can set 
\[
\eta_\tau:=\frac{\eta}{\eta(\tau)},\  \Phi_\tau=\Phi-\Phi(\tau)\otimes \eta_\tau,\ g_\tau=\eta_\tau^2\oplus \frac{1}{\eta(\tau)}g|_F.
\]
In particular, for any compact Sasakian $(X,g,\xi)$, there is another Sasakian structure, corresponding to a quasi-regular Reeb vector $\tau\in \mathfrak t_+$, on $X$ with biholomorphic Kähler cone $\mathcal C^*$ which is quasi-regular. By \cite[Theorem 1.3]{Ornea2007}, there is a projective orbifold $Y_\tau$ with an ample orbifold line bundle $L_\tau$ such that $\mathcal C:=\mathcal C\cup \{o\}$ is $\mathbb C^*$-equivariantly isomorphic to the affine variety
\[
\Spec\left(\bigoplus_{k\geq 0}H^0(Y_\tau,L_\tau^k)\right),
\]
and the induced map $\pi: X\to Y_\tau$ has the property that 
\begin{equation}\label{pullbeck of Chern class to the link}
    \pi^*c_1(L_\tau)=c_\tau \cdot [d\eta_\tau]\in H^{1,1}_B(X,\tau)(\simeq H^{1,1}_{orb}(Y_\tau)).
\end{equation}
for some positive constant $c_\tau$. Let $\mu: Y_\tau \to \bar Y_\tau$ be the coarse moduli space morphism. By the description of Seifert $\mathbb C^*$-bundles by \cite[Definition 5, Theorem 7]{KollarSeifert}, there exists an ample $\mathbb Q$-Cartier divisor $D$ such that  
\[
\mu_*L_\tau^k\simeq \mathcal O_{\bar Y_\tau}(\lfloor kD\rfloor), \qquad (\forall k\in \mathbb Z).
\]
Thus, we have 
\[
H^0(Y_\tau, L_\tau^k)\simeq H^0(\bar Y_\tau,\mu_*L_\tau^k)=H^0(\bar Y_\tau,\mathcal O_{\bar Y_\tau}(\lfloor kD\rfloor)).
\]
We recall the following standard fact, whose short proof is included for completeness:
\begin{lem}
Let $Y$ be an integral normal projective variety and $D$ be an ample $\mathbb Q$-Cartier divisor. Then the section ring
\[
R(Y,D)=\bigoplus_{m\geq 0}H^0(Y,\mathcal O_Y(\lfloor mD\rfloor))t^m
\]
is a finitely generated normal ring.
\end{lem}
\begin{proof}
The finite generation follows by passing to a Veronese subalgebra. 
For each prime divisor $E$ on $Y$, define a Gauss valuation $v_E$ on $K(Y)(t)=K(R(Y,D))$ by the extension of
\[
\tilde v_E\left(\sum_m f_mt^m\right)=\min_{m}(v_E(f_m)+m\ord_E(D)),
\]
where $v_E$ is the valuation on $K(Y)$ corresponding to $E$. Then we have 
\[
\tilde v_E(f_mt^m)\geq 0\iff v_E(f_m)\geq -m\ord_E(D)\iff f_m\in \mathcal O_Y(\lfloor mD\rfloor)\text{ near the generic point of $E$}.
\]
The last equivalence follows from the fact that $v_E$ is $\mathbb Z$-valued. Thus we have 
\[
R(Y,D)=K(Y)[t]\cap \bigcap_{E}\mathcal O_{\tilde v_E}
\]
and hence it is normal.
\end{proof}
Since $\bar Y_\tau$ is normal, we conclude that $\mathcal C=\Spec(\bigoplus_{k\geq 0}H^0(Y_\tau,L_\tau^k))$ is normal.
Since $T$ is commutative, its action on $X$ commutes with the $T_\tau$-action and hence descends to an action on $Y_\tau$. Moreover, $L_\tau$ inherits a $T$-linearization. Therefore, for each $k$, we obtain a representation
\[
T\to \GL(H^0(Y_\tau,L_\tau^k)).
\] 
Thus we have a weight decomposition 
\[
\mathcal C=\Spec \left(\bigoplus_{\chi\in \Hom(T_\mathbb C,\mathbb C^*)}R_\chi\right).
\]
There is a natural identification of $\Hom(\frac{1}{2\pi}\Lambda,\mathbb Z) \simeq \Hom(T_\mathbb C,\mathbb C^*)=:\mathbb X^*(T_\mathbb C)$. Define $\omega\subset \mathbb X^*(T_\mathbb C)_\mathbb R$ as the cone generated by $\chi\in \mathbb X^*(T_\mathbb C)$ such that $R_\chi \neq 0$. Then we can show that $\omega$ is a closed convex rational polyhedral cone since $R=H^0(\mathcal C,\mathcal O_\mathcal C)$ is finitely generated. We denote 
\[
\sigma=\omega^\vee:=\{\tau\in \mathfrak t\mid \langle \chi,\tau\rangle\geq 0, \forall \chi\in \omega\}.
\]
By the standard identification of the Reeb cone with the dual of the weight cone (see: \cite[Lemma 2.10, Lemma 2.13]{APOSTOLOV2021107969}), we have
\[
\mathfrak t_+={\rm Int}(\sigma).
\]

$c_\tau$ in equation (\ref{pullbeck of Chern class to the link}) can be computed as follows: Let $P_\tau$ be the minimal period of the Reeb flow generated by $\tau$, i.e.
\[
P_\tau=\int_{\text{a general $\tau$-orbit}}\eta_\tau.
\]
Let $S^1\times X\to X$ be the $S^1$-action generated by $\frac{P_\tau}{2\pi}\tau$. Then the induced connection $\sqrt{-1}A_\tau=\frac{2\pi\sqrt{-1}}{P_\tau}\eta_\tau$ on $S^1$-orbifold bundle $X\to Y$ coincides with that of $L^{-1}$. Thus we have 
\[
\pi^*c_1(L)=\frac{1}{P_\tau}[d\eta_\tau].
\]
Then we have $c_\tau=P_\tau^{-1}$ and thus
\[
\vol(X,\tau)=\int_X(d\eta_\tau)^n\wedge \eta_\tau=P_\tau c_\tau^{-n}c_1(L_\tau)^n=P_\tau^{n+1}c_1(L_\tau)^n.
\]

Let $E$ be a $T$-linearized transverse holomorphic vector bundle on Sasakian manifold $(X,g,\xi)$. 
As the Reeb vector varies in $\mathfrak t_+$, $E$ acquires an induced transverse holomorphic structure with respect to each new Reeb vector.
\begin{lem}\label{moving Reeb vector}
    Let $E$ be a $T$-linearized transverse holomorphic vector bundle on a compact Sasakian manifold $(X,g,\xi)$. Let $\tau\in \mathfrak t_+$. Define a partial connection $D''_\tau$ in the direction of $F^{0,1}\oplus \langle \tau\rangle$ as
    \[
    D''_\tau=D''_\xi|_{F^{0,1}}+\mathcal L_\tau\otimes \eta_\tau.
    \]
    Then $D''_\tau$ defines a transverse holomorphic structure of $E$ on $(X,g_\tau,\tau)$ and $T$ acts transverse holomorphically on $(E,D''_\tau)$.
\end{lem}
\begin{proof}
    Note that $(X,g_\tau,\tau)$ has the same underlying CR structure. Let $\bar \partial_{E,b}:=D''_\xi|_{F^{0,1}}$. Since $T$ acts transversely holomorphically on $E$, 
    \begin{equation}\label{transverse holomorphic action}
        g(\bar\partial_{E,b}(s)(v))=\bar\partial_{E,b}(gs)(gv)
    \end{equation} 
    for any $g\in T$, any local section $s$ of $E$, and any local vector field $v\in F^{0,1}$. Differentiating (\ref{transverse holomorphic action}) in the direction of $\tau$, we obtain
    \[
    \mathcal L_\tau \bar\partial_{E,b,v}(s)=\bar\partial_{E,b,v}(\mathcal L_\tau s)+\bar\partial_{E,b,[\tau,v]}(s).
    \]
    On the other hand, we have 
    \[
    R(D''_\tau)(u,v)s=D''_{\tau,u}D''_{\tau,v}s-D''_{\tau,v}D''_{\tau,u}s-D''_{\tau,[u,v]}s.
    \]
    If $u,v\in F^{0,1}$, we have $[u,v]\in F^{0,1}$ by the integrability of the CR structure, hence $R(D''_\tau)(u,v)=0$. If $v\in F^{0,1}$, we have 
    \[
    R(D''_\tau)(\tau,v)s=\mathcal L_\tau(\bar\partial_{E,b,v}(s))-\bar \partial_{E,b,v}(\mathcal L_\tau s)-\bar\partial_{E,b,[\tau,v]}s=0.
    \]
    Thus we have $R(D''_\tau)=0$. Since $T$ is commutative, $L_\tau(gs)=gL_\tau s$ for any $g\in T$ and local section $s$. Furthermore, the holomorphicity of $T$ action on $(E,D''_\xi)$ shows that $D''_\xi|_{F^{0,1}}(gs)=gD''_\xi|_{F^{0,1}}(s)$ for any $g$ and $s$. Thus the statement follows.
\end{proof}

Pulling back $E$ to the open cone $\mathcal C^*=X\times \mathbb R_{>0}$, we have a $T_\mathbb C$-linearized holomorphic vector bundle $\tilde E$. 
\begin{lem}
    There is a $T_\mathbb C$-linearized algebraic reflexive sheaf $\mathcal E$ on the closed cone $\mathcal C$ such that $(\mathcal E|_{C^*})^{{\rm an}}=\tilde E$. 
\end{lem}
\begin{proof}
    Let $\tau \in \Lambda \cap \mathfrak t_+$. Then $E_\tau$ is an $S^1$-linearized transverse holomorphic vector bundle on $(X,\tau)$, and hence descends to a holomorphic orbifold vector bundle $\bar E$ on $Y_\tau=[X/S^1]$. Since $Y_\tau$ is projective DM-stack, $\bar E$ is algebraic by orbifold GAGA \cite[Theorem 2.4]{Hall_2014}. Since $T$ is commutative and contains $T_\tau$, $T$-linearization of $E$ descends to the $T$-action on $\bar E$, which induces a representation 
    \[
    T\to \GL(H^0(Y_\tau,\bar E\otimes L_\tau^k))
    \]
    for any $k\in \mathbb Z$ since each $g\in T$ acts holomorphically. Thus we have a weight decomposition
    \[
    H^0(Y_\tau,\bar E\otimes L_\tau^k)=\bigoplus_{\chi}H^0(Y_\tau,\bar E\otimes L_\tau^k)_\chi.
    \]
    As we noted above, $L_\tau$ also inherits $T$-linearization. By these constructions, multiplication map
    \[
    H^0(Y_\tau,L_\tau^l)\otimes H^0(Y_\tau,\bar E\otimes L_\tau^k)\to H^0(Y_\tau,\bar E\otimes L_\tau^{k+l})
    \]
    is $T$-equivariant, and hence we get the $\mathbb X^*(T_\mathbb C)$-graded $R=\Gamma(\mathcal C,\mathcal O_\mathcal C)$-module
    \[
    M:=\bigoplus_{\chi\in \mathbb X^*(T_\mathbb C)}M_\chi, \qquad M_\chi=\bigoplus_k H^0(Y_\tau,\bar E\otimes L_\tau^k)_\chi.
    \]
    It is finitely generated since $L_\tau$ is ample. Since $\bar E$ and the natural projection $\pi: \mathcal C^*=\mathcal C\setminus \{0\}=\Spec(\bigoplus H^0(Y,L^k))\setminus \{0\}\to Y_\tau$ are algebraic, $\tilde E=\pi^*\bar E$ is algebraic. Since the codimension of $\{0\}$ in $\mathcal C$ is $\geq 2$ and $\mathcal C$ is normal, the extension $\mathcal E:=j_*\tilde E$ is a reflexive sheaf on $\mathcal C$. From the weight decomposition above, the $T_{\mathbb C}$-linearization on $\bar E$ extends to that of $\tilde E$. By construction, we have 
    \[
    H^0(\mathcal C,\mathcal E)=M
    \]
    By construction, the analytification of $\mathcal E$ restricts to $\tilde E$ on $\mathcal C^*$. This proves the statement.
\end{proof}

In \cite{hada2026optimalalgebraictangentcone}, we introduced a notion of stability for $\mathbb R$-graded modules over $\mathbb R$-graded domains as follows:
\begin{propdef}(=\cite[Theorem 2.2, Definition 2.13]{hada2026optimalalgebraictangentcone})\label{asymptotic RR and definition of graded module stability}
    Let $R=\bigoplus_{\gamma\in \mathbb R_{\geq 0}}R_\gamma$ be a finitely generated graded domain over $R_0=\mathbb C$, and $M=\bigoplus_{\lambda \in \mathbb R}M_\lambda$ be a finitely generated $\R$-graded $R$-module. Then we have
    \[
    \sum_{\lambda \leq x}\dim M_\lambda=\frac{a_d(M)}{d!}x^d+\frac{a_{d-1}(M)}{(d-1)!}x^{d-1}+R(x),
    \]
    with 
    \[
    R(x)=O(x^{d-1}), \ \int_0^TR(x)dx=o(T^d),
    \]
    where $d=\dim(M)$. Furthermore, if $d=\dim(R)$, we have 
    \[
    a_d(M)=\rank (M)\cdot a_d(R).
    \]
    Define $\deg(M)=a_{d-1}(M)-\rank(M)\cdot a_{d-1}(R)$, and 
    \[
    \mu_{\gr}(M):=\frac{\deg(M)}{a_d(M)}
    \]
    for a torsion-free graded module $M$. $M$ is stable (resp. semistable) if $\mu_{\gr}(N)<\mu_{\gr}(M)$ (resp. $\mu_{\gr}(N)\leq \mu_{\gr}(M)$) for any saturated graded submodule $N\subset M$ with $0<\rank N<\rank M$. 
\end{propdef}

Let $R$ be a $\mathbb X^*(T_\mathbb C)$-graded ring and $M$ be a $\mathbb X^*(T_\mathbb C)$-graded $R$-module. For $\xi\in \mathfrak t$, coarsening the $\mathbb X^*(T_\mathbb C)$-grading along 
\[
\chi\mapsto \langle \chi,\xi\rangle
\]
gives an $\mathbb R$-graded ring $R_\xi$ and an $\mathbb R$-graded $\mathbb R_\xi$-module $M_\xi$.

Our goal is to prove the following theorem:
\begin{Thm}
    Let $E$ be a $T$-linearized transverse holomorphic vector bundle on a compact Sasakian manifold $(X,g,\xi)$. Then we have
    \[
    \mu_{\Sas}(E)=\frac{1}{2\pi}\mu_{\gr}(M_\xi),
    \]
    where $M$ is defined above. Furthermore, $E$ is (semi)stable if and only if $M_\xi$ is (semi)stable.
\end{Thm}
\section{Equivalence of two stability conditions}
\subsection{Slope identity}
Let $(X,g,\xi)$ be a quasi-regular compact Sasakian manifold of dimension $2n+1$ for $n\geq 1$. Let $E$ be a transverse holomorphic vector bundle with $T=T_\xi\simeq S^1$-action. 
Then $E$ descends to an orbifold vector bundle $\bar E$ on $Y=[X/T_\xi]$, and the corresponding graded module $M$ is given by
\[
M=\bigoplus_{\chi\in \mathbb X^*(T_\mathbb C)}H^0(Y,\bar E\otimes L^{m(\chi)}),
\]
where $L$ is the ample orbifold line bundle with the following property:
\[
\mathcal C=\Spec\left(\bigoplus_{\chi\in \mathbb X^*(T_\mathbb C)}H^0(Y,L^{m(\chi)})\right),
\]
and $m(\chi)$ is defined as 
\[
m(\chi)=\frac{1}{2\pi}\langle P_\xi \xi,\chi\rangle\in \mathbb Z.
\]
It is natural to regard $R$ and $M$ as $\{\langle \xi,\chi\rangle\mid \chi\in \mathbb X^*(T_\mathbb C)\}=\frac{2\pi}{P_\xi}\mathbb Z$-graded. By the orbifold Riemann-Roch \cite{Kaw79}, we have 
\[
\deg(M)=\left(\frac{P_\xi}{2\pi}\right)^n\int_Y c_1(\bar E)\wedge c_1(L)^{n-1}.
\]
Indeed, since $L$ is the orbifold line bundle associated with the Seifert $S^1$-bundle $X\to Y$, $g\neq1$ acts nontrivially on its fiber.

Fix a nontrivial element $g$ occurring on an inertia component $F$. Then $g$ corresponds to some root of unity $\zeta\neq 1$ and thus
\[
{\rm ch}_g(L^k)=\zeta^ke^{kc_1(L|_F)}.
\]
Thus, its contribution is
\[
S_N:=\sum_{0\leq k\leq N}\zeta^kP(k)=a\zeta^NN^n+O(N^{n-1}),
\]
where $P$ is a polynomial of degree at most $n$ and $a$ is a constant depending on the leading coefficient of $P$. Thus we have
\[
\lim_{N\to \infty}\frac{\sum_{0\leq k\leq N}S_k}{N^{n+1}}=0.
\]
Thus we have
\begin{align*}
    \deg(M)=&\left(\frac{P_\xi}{2\pi}\right)^n\int_Yc_1(\bar E)\wedge c_1(L)^{n-1}\\
    =&\left(\frac{1}{2\pi}\right)^n\int_Xc_1(E)\wedge (d\eta)^{n-1}\wedge \eta\\
    =&\left(\frac{1}{2\pi}\right)^n\vol(X)\deg(E)\\
    =&2\pi a_{n+1}(R_\xi)\deg(E).
\end{align*}
Since $\rank M=\rank E$, we have the following lemma:
\begin{lem}\label{slope identity in quasi-regular case}
    Let $(X,g,\xi)$ be a quasi-regular compact Sasakian manifold. For any transverse holomorphic vector bundle $E$ on $(X,g,\xi)$, we have the following equality
    \[
    \mu_{\gr}(M)=2\pi\cdot \mu_{\Sas}(E),
    \]
    where $M$ is the associated $\mathbb R$-graded $R$-module.\qed
\end{lem}
When $\xi$ is irregular, a $T$-linearized transverse holomorphic vector bundle $E$ on $X$ corresponds to an $\mathbb X^*(T_\mathbb C)$-graded module
\[
M=\bigoplus_{\chi\in \mathbb X^*(T_\mathbb C)}M_\chi
\]
over $R=\bigoplus_{\chi\in \mathbb X^*(T_\mathbb C)}R_\chi$. For $\tau\in \mathfrak t_+$, define the $\mathbb R$-graded algebra $R_\tau$ and the $\mathbb R$-graded $R_\tau$-module $M_\tau$ by
\[
R_\tau=\bigoplus_{\lambda\in \mathbb R}\left(\sum_{\langle \tau,\chi\rangle=\lambda}R_\chi\right)
\]
and 
\[
M_\tau=\bigoplus_{\lambda\in \mathbb R}\left(\sum_{\langle \tau,\chi\rangle=\lambda}M_\chi\right).
\]
Let $v_\tau: R\to \mathbb R$ be a valuation defined as
\[
v_\tau(f)=\min\{\langle \chi,\tau\rangle\mid f_\chi\neq 0\},
\]
where $f=\sum_{\chi} f_\chi$ is a weight decomposition. Then we have 
\[
\vol(R_\tau)=\vol(v_\tau)=\lim_{m\to \infty}\frac{\dim(R/\mathcal J_{\geq m}(v_\tau))}{m^{n+1}/(n+1)!}
\]
For the arguments in the algebraic setting, we slightly generalize the setup. Let $\mathcal C$ be a normal affine variety with a good $T_\mathbb C$-action in the sense of \cite[Definition 2.8]{li2021algebraicity}; namely, an effective $T_\mathbb C$-action with a closed point $o\in \mathcal C$ such that 
\begin{itemize}
    \item $o$ is fixed by $T_\mathbb C$-action. 
    \item For any $x\in \mathcal C$, $o$ is in the closure of $T_\mathbb C\cdot x$.
\end{itemize}
Then, we have the weight decomposition of $R=\Gamma(\mathcal C,\mathcal O_{\mathcal C})$:
\[
R=\bigoplus_{\chi\in \mathbb X^*(T_\mathbb C)}R_\chi,
\]
The Reeb cone associated with a good $T_{\mathbb C}$-action is defined in \cite{li2021algebraicity} in the same way. Note that \cite{li2021algebraicity} does not impose normality or integrality in the definition of a good $T_{\mathbb C}$-action, while in this paper we work only with integral normal affine varieties.

In the algebraic setting, we consider torsion-free $\mathbb X^*(T_\mathbb C)$-graded modules (or their reflexive hulls) instead of $T$-linearized transverse holomorphic vector bundles on compact Sasakian manifolds.

We mainly formulate our results in the Sasakian setting. However, several statements and arguments below remain valid in the purely algebraic setting described above. We mark such statements by $\diamond$. When a marked statement is formulated in Sasakian terms, its algebraic counterpart is obtained by replacing a $T$-linearized transverse holomorphic vector bundle by a torsion-free $\mathbb X^*(T_{\mathbb C})$-graded module. Note that in the algebraic setting slope of torsion-free sheaf is the same as that of reflexive hull.

\begin{alglemma}
    In both the algebraic and Sasakian settings, the functions $\tau\in \mathfrak t_+\mapsto \vol(R_\tau)$ and $\tau\in \mathfrak t_+\mapsto \deg(M_\tau)$ are smooth.
\end{alglemma}
\begin{proof}
    \cite[Lemma 2.11]{li2021algebraicity} proves the smoothness of the volume function. The same argument applies here; for completeness, we include the details. 
    
    Let $F(\tau,t)$ be the index character, i.e.
    \[
    F(\tau,t)=\sum_{\chi}e^{-t\langle \chi,\tau\rangle}\dim M_\chi=\int_0^\infty te^{-tx}f_\tau(x)dx+O(1),
    \]
    where 
    \[
    f_\tau(x)=\sum_{\langle \chi,\tau\rangle\leq x}\dim M_\chi.
    \]
    By \cite[Theorem 2.2]{hada2026optimalalgebraictangentcone}, we have 
    \[
    f_\tau(x)=\frac{a_{n+1}(M_\tau)}{(n+1)!}x^{n+1}+\frac{a_n(M_\tau)}{n!}x^n+R_\tau(x)
    \]
    with $R_\tau(x)=O(x^n)$ and 
    \[
    \lim_{T\to \infty} \frac 1 {T^{n+1}} \int_0^T R_\tau(x)dx=0.
    \]
    Thus we have 
    \[
    t^{n+1}F(\tau,t)=a_{n+1}(M_\tau)+ta_n(M_\tau)+t^{n+2}\int_0^\infty e^{-tx}R_\tau(x)dx+O(t^{n+1}).
    \]
    Let 
    \[
    F_\tau(T):=\int_0^T R_\tau(x)dx.
    \] 
    Then we have 
    \begin{align*}
        t^{n+1}\int_0^\infty e^{-tx}R_\tau(x)dx=&t^{n+2}\int_0^\infty e^{-tx}F_\tau(x)dx\\
        =&t^{n+1}\int_0^\infty e^{-u}F_\tau\left(\frac{u}{t}\right)du.
    \end{align*}
    Since $F_\tau(x)=o(x^{n+1})$, $t^{n+1}F_\tau(u/t)\to 0$ ($t\to 0$). Furthermore, there is some $C_\tau>0$ such that $|F_\tau(x)|\leq C_\tau (1+x^{n+1})$ ($x\geq 0$) by $R_\tau(x)=O(x^n)$. By the dominated convergence theorem, this integral tends to $0$. Thus we have the following asymptotic expansion of $F(\tau,t)$:
    \[
    F(\tau,t)=\frac{a_{n+1}(M_\tau)}{t^{n+1}}+\frac{a_n(M_\tau)}{t^n}+o(t^{-n}).
    \]
    By \cite[Proposition 4.9]{CS18}, the multigraded Hilbert series of $M$ is rational. After the specialization $z^\alpha=e^{-t\langle  \alpha,\tau\rangle}$, this yields a rational expression for the index character $F_M(\tau,t)$. The argument in the proof of \cite[Theorem 4.10]{CS18} then shows that its Laurent coefficients depend smoothly on $\tau\in\mathfrak t_+$. 
    Hence $a_{n+1}(M_\tau)$ and $a_n(M_\tau)$ depend smoothly on $\tau$. Applying the same argument to $M$ and $R$, we conclude that $\vol(R_\tau)$ and $\deg(M_\tau)$ are smooth as well. 
\end{proof}
Since $\vol(R_\tau)>0$ for any $\tau\in \mathfrak t_+$, we have the following corollary:
\begin{algcorollary}\label{smoothness of graded slope}
    Both in algebraic and Sasakian setting, the function $\tau\in \mathfrak t_+\mapsto \mu_{\gr}(M_\tau)$ is smooth.\qed
\end{algcorollary}

This suggests that the Reeb cone should admit a decomposition according to Harder--Narasimhan filtrations. We return to this question in Section 4.
By construction, the $\mathbb R$-graded module associated with $E_\tau:=(E,D''_\tau)$, where $D''_\tau$ is constructed in Lemma \ref{moving Reeb vector}, is $M_\tau$, since neither $D''_\tau|_{F^{0,1}}$ nor $T$-action changes. Fix an arbitrary $T$-invariant Hermitian metric $h$ on $E$. Then the Chern connection $D_\tau$ with respect to $(D''_\tau,h)$ depends smoothly on $\tau\in \mathfrak t_+$ since $D''_\tau$ depends smoothly on $\tau\in \mathfrak t_+$, and hence its Chern form $c_1(E_\tau,h)$ also depends smoothly on $\tau$. Moreover, $\eta_\tau$, and hence $\vol(X,\tau)$, depends smoothly on $\tau$. As a consequence, we have the following theorem:
\begin{Thm}\label{equality of slopes}
    Let $(X,g,\xi)$ be a compact Sasakian manifold, and $E_\xi=(E,D''_\xi)$ be a $T$-linearized transverse holomorphic vector bundle on $(X,g,\xi)$. Then for any $\tau\in \mathfrak{t}_+$, we have
    \[
    \mu_{\Sas}(E_\tau)=\frac{1}{2\pi}\mu_{\gr}(M_\tau),
    \]
    where $M_\tau$ is an $\mathbb R$-graded $R_\tau$-module associated with $E_\tau$.
\end{Thm}
\begin{proof}
    Lemma \ref{slope identity in quasi-regular case} shows that for any $\tau\in \mathfrak t_{\mathbb Q,+}$, 
    \[
    \mu_{\Sas}(E_\tau)=\frac{1}{2\pi}\mu_{\gr}(M_\tau),
    \]
    and both sides depend smoothly on $\tau\in \mathfrak t_+$. Since $\mathfrak t_{\mathbb Q,+}\subset \mathfrak{t}_+$ is dense, we have the required consequence.
\end{proof}
\subsection{Subsheaf-submodule correspondence}
In this section, we prove the following theorem:
\begin{Thm}\label{subsheaf-submodule correspondence}
    Let $(X,g,\xi)$ be a compact Sasakian manifold, and $(E,D''_\xi)$ a $T$-linearized transverse holomorphic vector bundle on $(X,g,\xi)$, and $M$ the corresponding $\mathbb X^*(T_\mathbb C)$-graded module. Then there is a canonical correspondence between saturated transverse coherent subsheaves of $(E,D''_\xi)$ and saturated $\mathbb X^*(T_\mathbb C)$-graded submodules of $M$.
\end{Thm}

\begin{Rem}
    Since the Reeb flow is dense in $T$, the homomorphism 
    \[
    \mathbb X^*(T_\mathbb C)\to \mathbb R; \qquad \chi\mapsto \langle \chi,\xi\rangle
    \]
    is injective. Hence, this correspondence may be regarded as the correspondence between saturated transverse coherent subsheaves of $(E,D''_\xi)$ and saturated $\mathbb R$-graded submodules of $M_\xi$.
\end{Rem}
\begin{proof}
    Let $\mathcal F$ be a saturated transverse coherent subsheaf of $E$. Then $\mathcal F$ is $T$-invariant by Lemma \ref{saturated subsheaves of $T$-linearized bundle are $T$-invariant}. Take $\tau\in \mathfrak t_{\mathbb Q,+}$, and let $p_\tau:X\to Y_\tau$ be the canonical morphism. Descending the inclusion $\mathcal F\subset E$ to $Y_\tau=[X/T_\tau]$, we obtain an inclusion $\bar{\mathcal F}\subset \bar E$. Then we have natural inclusions $\bar{\mathcal F}\otimes L_\tau^k\subset \bar E\otimes L_\tau^k$, and hence 
    \[
    \bigoplus_{k\in \mathbb Z}H^0(Y_\tau, \bar {\mathcal F}\otimes L_\tau^k)\subset \bigoplus_{k\in \mathbb Z}H^0(Y_\tau,\bar E\otimes L_\tau^k).
    \]
    Hence, the same construction as in Section 2.2 induces a $\mathbb X^*(T_\mathbb C)$-graded submodule $N_\mathcal F\subset M$. Assume that $N_\mathcal F\subset M$ is not saturated. Then $\bigoplus_k H^0(Y_\tau,(\bar E/\bar{\mathcal F})\otimes L_\tau^k)$ has torsion since there is natural inclusion
    \[
    \frac{H^0(\bar E\otimes L_\tau^k)}{H^0(\bar{\mathcal F}\otimes L_\tau^k)}\subset H^0((\bar E/\bar{\mathcal F})\otimes L_\tau^k).
    \]
    Thus there is some $0\neq s\in H^0(Y_\tau,(\bar E/\bar{\mathcal F})\otimes L_\tau^k)$ and $0\neq f\in H^0(Y_\tau,L_\tau^l)$ such that $fs=0$. For any $x\in X$, $(p_\tau^*L_\tau^l)_x\simeq \mathcal O_{X,x}$ and thus we have $\tilde f_x\cdot \tilde s_x=0$, where $\tilde f_x$ is the corresponding holomorphic germ to $f_x$ and $\tilde s_x$ is the image of $p_\tau^*s_x$ by the natural isomorphism $p_\tau^*(\bar E/\bar {\mathcal F})\to E/\mathcal F$. Choose $x\in X$ such that $f_{p_\tau(x)}\neq 0$ and $s_{p_\tau(x)}\neq 0$. Then both $\tilde f_x$ and $\tilde s_x$ are nonzero. This contradicts the saturatedness of $\mathcal F\subset E$. The analytification of the restriction of the corresponding sheaf $(N_\mathcal F)^\sim$ on $\mathcal C$ to $\mathcal C^*$ is $T_\mathbb C$-isomorphic to $\pi^*\mathcal F$, where $\pi: \mathcal C^*\to X$ is natural projection. Thus this construction doesn't depend on the choice of $\tau$.
    
    On the other hand, for any saturated graded submodule $N\subset M$, we associate a saturated transverse coherent subsheaf $\mathcal F_N\subset E$ as follows: Take $\tau\in \mathfrak t_{\mathbb Q,+}$. Then $N_\tau\subset M_\tau\simeq \bigoplus_{k}H^0(Y_\tau,\bar E\otimes L_\tau^k)$ is finitely generated submodule. Take homogeneous generators $s_i$ of $N$, and hence of $N_\tau$, and regard $s_i$ as a section of $\bar E\otimes L_\tau^{k_i}$. Define a sheaf $\mathcal F_N$ as the $p_\tau^*\bar {\mathcal F}_N$. Then $\bar {\mathcal F}_N$ is defined by the following:
    \[
    \bar{\mathcal F}_N:={\rm Im}\left((s_i): \bigoplus_{i} L_\tau^{-k_i}\to \bar E\right).
    \]
    This map $(s_i)$ is naturally $T$-equivariant. Thus $\mathcal F_N$ is $T$-equivariant transverse coherent subsheaf of $E$. If $\mathcal F_N\subset E$ is not saturated, $\pi^*(E/F_N)=((M/N)^\sim|_{\mathcal C^*})^{{\rm an}}$ has torsion, hence $N\subset M$ is not saturated. Thus $\mathcal F_N\subset E$ is saturated. This construction can be rephrased as follows: for saturated $N\subset M$, restrict $(N^\sim)^{{\rm an}}\subset \mathcal E^{{\rm an}}$ to the link $X\simeq L_\xi=\{r=1\}\subset \mathcal C$. Thus this construction doesn't depend on the choice of $\tau$. 
    
    It is easy to see that these constructions are inverse to each other since saturated subsheaf of $M^\sim|_{\mathcal C^*}$ extends uniquely to the saturated subsheaf of $M^\sim$ and its global sections coincide with those on $\mathcal C^*$ by reflexivity and normality of $\mathcal C$. Thus we get the required statement.
\end{proof}

\subsection{Equivalence of stabilities}
In this section, we prove the following theorem:
\begin{Thm}\label{equivalence of stabilities}
    Let $(X,g,\xi)$ be a compact Sasakian manifold of dimension $\geq 3$. For any $T$-linearized transverse holomorphic vector bundle $(E,D''_\xi)$ on $(X,g,\xi)$, the following conditions are equivalent:
    \begin{enumerate}
        \item $(E,D''_\xi)$ is (semi)stable in the sense of Definition \ref{BS-stability}.
        \item The associated $\mathbb R$-graded $R_\xi$-module $M_\xi$ is (semi)stable in the meaning of Definition \ref{asymptotic RR and definition of graded module stability}.
    \end{enumerate}
\end{Thm}

To prove the theorem, we first establish the following lemma.
\begin{alglemma}\label{determinant degree equals the initial degree}
    Let $R=\bigoplus_{\gamma\in \mathbb R_{\geq 0}}R_\gamma$ be a normal $\mathbb R_{\geq 0}$-graded domain of finite type over $R_0=\mathbb C$ of dimension $n+1$. For any finitely generated torsion-free $\mathbb R$-graded $R$-module $M$, we have 
    \[
    \deg(M)=\deg(\det(M)),
    \]
    where $\det(M)=(\wedge^rM)^{**}$ is the determinant module.
\end{alglemma}
\begin{proof}
    By \cite[Lemma 2.3]{hada2026optimalalgebraictangentcone}, we can construct an exact sequence
    \[
    0 \to \bigoplus_{i=1}^rR(-\lambda_i)\overset{\phi}\to M\to Q\to 0
    \]
    of $\mathbb R$-graded $R$-modules with torsion module $Q$. Then we have 
    \[
    \deg(M)=-\sum_{i=1}^r\lambda_i a_{n+1}(R)+\sum_{P}\length(H^0_P(Q_P))a_n(R/P),
    \]
    where $P$ in the second sum runs over homogeneous prime ideals of height $1$. 

    Taking top exterior power, we have an injection
    \[
    R\left(-\sum_i \lambda_i\right)\to \wedge^rM.
    \]
    Let $Q'$ be the cokernel of this map. Take a homogeneous prime $P$ of height $1$. Then $R_P$ is DVR by normality. Let $\pi_P$ be a uniformizer of $R_P$. Thus $\phi_P$ is a morphism of free modules and hence we have $m_1\dots,m_r$ such that $\phi_P={\rm diag}(\pi_P^{m_1},\dotsm\pi_P^{m_r})$ with respect to some appropriate basis. 
    \[
    Q_P\simeq \bigoplus_{i=1}^rR_P/\pi_P^{m_i}
    \]
    as $R_P$-module. Thus $\length Q_P=\sum_{i=1}^r m_i$. Taking the top exterior product yields a morphism $\pi_P^{\sum_i m_i}: R_P\to R_P$ and hence
    \[
    Q'_P\simeq R_P/\pi_P^{m_1+\dots +m_r},
    \]
    Thus we have 
    \[
    \length Q_P=\length Q'_P.
    \]
    This shows that $\deg(M)=\deg(\wedge^rM)$. Since $\det(M)=(\wedge^r M)^{**}$ and $\wedge^r M$ are isomorphic in codimension $1$, $\deg(\det(M))=\deg(\wedge^r M)$. Therefore, we have the required consequence.
\end{proof}

\noindent
{\it Proof of Theorem \ref{equivalence of stabilities}.}

Let $\mathcal F\subset E$ be a saturated subsheaf of rank $0<s<\rank E=:r$. Then we have a saturated subsheaf $\det(\mathcal F)\subset \wedge^s E$. Thus we have an associated $\mathbb R$-graded $R$-module $L$. By looking at the weight decomposition, we have $L=\det(N_\mathcal F)$. Therefore, we have 
\[
\mu_{\Sas}(\mathcal F)=\frac{1}{2\pi}\mu_{\gr}(N_\mathcal F)
\]
by Lemma \ref{determinant degree equals the initial degree} and Theorem \ref{equality of slopes} since $\det(\mathcal F)$ is $T$-linearized transverse holomorphic line bundle. Together with Theorem \ref{subsheaf-submodule correspondence}, this shows that $E$ is (semi)stable if and only if $M_\xi$ is (semi)stable.\qed

\section{Wall structures for stability}
Let $(X,g,\xi)$ be a compact Sasakian manifold of dimension $2n+1$ ($n\geq 1$). Let $(E,D''_\xi=\bar\partial_{E,b}+\mathcal L_\xi\otimes \eta)$ be a $T$-linearized transverse holomorphic vector bundle on $(X,g,\xi)$. Then there is an associated $\mathbb X^*(T_\mathbb C)$-graded module $M$ over $\mathbb X^*(T_\mathbb C)$-graded algebra $R$ of finite type over $R_0=\mathbb C$. In this section, we prove the openness of stability and discuss the structure of stable cones.

\subsection{Openness of stability: a differential-geometric proof}
In this subsection, we investigate the following openness result:
\begin{Thm}[Openness of stability I]\label{Openness of Stability I}
    Let $E$ be a $T$-linearized transverse holomorphic vector bundle on $(X,g,\xi)$. If $E_{\tau_0}$ is stable for some $\tau_0\in \mathfrak t_+$, there is an open neighbourhood $U\subset \mathfrak t_+$ such that $E_\tau$ is a stable transverse holomorphic vector bundle on $(X,g_\tau,\tau)$ for any $\tau\in U$. 
\end{Thm}
\begin{Rem}
A related openness result was recently obtained by Delloque \cite[Theorem 3.6]{Delloque25} in the setting of slope stability of torsion-free sheaves on balanced manifolds with respect to balanced classes. More precisely, for a fixed torsion-free coherent sheaf on a compact balanced manifold, he showed that the set of balanced classes with respect to which the sheaf is stable is an open convex cone, locally defined by finitely many linear inequalities.

Our setting is different in that the parameter is the Reeb vector rather than a balanced class. The underlying CR structure is fixed, and so is the complex structure of the associated affine cone, while the Reeb foliation and the transverse K\"ahler structure vary with $\tau\in\mathfrak t_+$. Algebraically, we keep the $\mathbb X^*(T_{\mathbb C})$-graded module fixed and vary the homomorphism
\[
X^*(T_{\mathbb C})_{\mathbb R}\to \mathbb R,
\qquad
\chi\mapsto \langle \chi,\tau\rangle .
\]
In contrast to the variation of balanced classes, the resulting slope functions are generally nonlinear in $\tau$. Thus the usual convexity argument by linear slope inequalities does not apply in general. In particular, when the affine cone is not $\mathbb Q$-factorial, the stable cone need not be convex and the walls may be nonlinear; see Example \ref{stable cone of some rank 2 bundle on T11}.
When the affine cone is $\mathbb Q$-factorial, however, the slope functions become linear and the stable cone is defined by finitely many linear inequalities; see Lemma \ref{liniarity of slopes when cone is Q-factorial} and Theorem \ref{finiteness of walls}.
\end{Rem}
To prove the theorem, we establish a first-derivative estimate for slope functions and a uniform spectral gap. We first prove a lemma that does not use differential geometry and is valid in both the algebraic and Sasakian settings.
\begin{alglemma}\label{First derivative estimate}
    Let $M$ be a finitely generated torsion-free $\mathbb X^*(T_\mathbb C)$-graded $R$-module. Fix a norm $\|\cdot \|$ on $\mathfrak t$ and its dual norm on $\mathfrak t^*=\mathbb X^*(T_\mathbb C)_\mathbb R$. For any $\xi\in \mathfrak t_+$, there is a constant $C_M>0$ depending only on $M$, and a positive constant $C'_\xi$, depending continuously on $\xi$ such that for any saturated graded submodule $N\subset M$ of $0<\rank N<\rank M$, the following inequality holds:
    \[
    \|d\mu(N_\tau)|_{\tau=\xi}\|\leq C_M+C'_\xi(C_M\|\xi\|-\mu(N_\xi)).
    \]
    Furthermore, $C'_\xi$ satisfies the following equality for any $\lambda>0$:
    \[
    C'_{\lambda\xi}=\lambda^{-1}C'_\xi
    \]
\end{alglemma}
\begin{proof}
    For each $0<s<\rank M$, take an embedding 
    \[
    \left(\bigwedge^sM\right)^{**}\hookrightarrow \bigoplus_{j=1}^{k_s}R(\beta_{s,j}),
    \]
    where $\beta_{s,j}\in \mathbb X^*(T_\mathbb C)$. Such an embedding can be constructed by taking the dual of the presentation of $(\wedge^sM)^*$. Let $N$ be a $\rank s$ saturated graded submodule of $M$. Then we have an inclusion
    \[
    \det(N)\hookrightarrow \left(\bigwedge^sM\right)^{**}\hookrightarrow \bigoplus_{j=1}^{k_s}R(\beta_{s,j}).
    \]
    There exists $j$ such that the induced morphism $\det(N)\to R(\beta_{s,j})$ is injective. Fix such a $j$ and define $Q$ by the exact sequence
    \[
    0\to \det(N)\to R(\beta_{s,j})\to Q\to 0
    \]
    Let $\ell_P(Q)=\length_P(Q_P)$ for height one homogeneous prime. Then we get
    \[
    \mu_\tau(N)=\frac{1}{s}\left(\langle \beta_{s,j},\tau\rangle-\sum_{P}\ell_P(Q)\Theta_P(\tau)\right),
    \]
    where 
    \[
    \Theta_P(\tau)=\frac{a_n((R/P)_\tau)}{a_{n+1}(R_\tau)}.
    \]
    We use the dual norm on $\mathfrak t^*=\mathbb X^*(T_\mathbb C)_\mathbb R$, also denoted by $\|\cdot \|$. For any $\xi\in \mathfrak t_+$, there exists $C_\xi>0$ such that 
    \[
    |\langle\chi,v\rangle|\leq C_\xi\|v\|\langle\chi,\xi\rangle
    \]
    for any $v\in \mathfrak t$ and $\chi\in \omega$. Indeed, $\langle\cdot,\xi\rangle: \mathbb X^*(T_\mathbb C)_{\mathbb R}\to \mathbb R$ is continuous, linear, and positive in the closed cone $\omega$. Thus we can define 
    \[
    C_\xi:=\max_{\substack{\chi\in \omega,\\ \|\chi\|=1}}\frac{1}{\langle \chi,\xi\rangle}.
    \]
    For a homogeneous height $1$ prime $P$ of $R$, define $F_P$ as its index character, i.e.
    \[
    F_P(\tau,t):=\sum_{\chi\in \mathfrak X^*(T_\mathbb C)}e^{-t\langle\chi,\tau\rangle}\dim(R/P)_\chi.
    \]
    Thus we have 
    \begin{align*}
        |D_vF_P(\xi,t)|\leq& \sum_{\chi}|t\langle\chi,v\rangle|e^{-t\langle\chi,\xi\rangle}\dim(R/P)_\chi\\
        \leq& C_\xi\|v\|\sum_\chi t\langle \chi,\xi\rangle e^{-t\langle\chi,\xi\rangle}\dim(R/P)_\chi\\
        =&C_\xi\|v\|\left(-t\frac{\partial}{\partial t}F_P(\xi,t)\right)
    \end{align*}
    for every $t>0$ and $v\in \mathfrak t$. Since $F_P(\xi,t)$ can be written as
    \[
    F_P(\xi,t)=\frac{a_n((R/P)_\xi)}{t^n}+O(t^{1-n}),
    \]
    we have 
    \[
    |D_va_n((R/P)_\xi)|\leq nC_\xi \|v\| a_n((R/P)_\xi).
    \]
    This means that 
    \[
    \|d\log a_n((R/P)_\xi)\|\leq nC_\xi
    \]
    Similarly, we have $\|d\log a_{n+1}(R_\xi)\|\leq (n+1)C_\xi$ and thus 
    \[
    \|d\log \Theta_P(\xi)\|\leq (2n+1)C_\xi=:C'_\xi.
    \]
    From this inequality, we have 
    \begin{align*}
        \|d\mu(N_\tau)|_{\tau=\xi}\|\leq& \frac{1}{s}\|\beta_{s,j}\|+\frac{C'_\xi}{s}\sum_{P}\ell_P(Q)\Theta_P(\xi)\\
        =&\frac{1}{s}\|\beta_{s,j}\|+C'_\xi\left(\frac{1}{s}\langle \beta_{s,j},\xi\rangle-\mu(N_\xi)\right)
    \end{align*}
    For any graded submodule $N'\subset M$ of rank $s'$, there exists $\beta_{s',j'}$ such that $\deg(N'_\xi)\leq \deg(R(\beta_{s',j'})_\xi)=a_{n+1}(R_\xi)\langle\beta_{s',j'},\xi\rangle$ and thus $\mu(N_\xi)\leq \frac{1}{s'}\langle \beta_{s',j'},\xi\rangle$. Set
    \[
    C_{M}:=\max_{s,j}\frac{1}{s}\|\beta_{s,j}\|.
    \]
    Then we have 
    \[
    \mu_{\max,\xi}(M):=\max\{\mu(N'_\xi)\mid N'\subset M\text{: $\mathbb X^*(T_\mathbb C)$-graded submodules}\}\leq C_M\|\xi\|.
    \]
    Thus we have 
    \[
    \|d\mu(N_\tau)|_{\tau=\xi}\|\leq C_M+C'_\xi\left(C_M\|\xi\|-\mu(N_\xi)\right),
    \]
    where $C'_\xi$ depends continuously. Note that $C'_{\lambda\xi}=\lambda^{-1}C'_\xi$, $\|\lambda\xi\|=\lambda\|\xi\|$, $\mu(N_{\lambda\xi})=\lambda\mu(N_\xi)$ for any $\lambda>0$. Thus the second term is invariant under rescaling.
\end{proof}

\begin{Thm}\label{Spectral Gap Theorem}
    Let $E_\xi$ be a transverse holomorphic vector bundle of rank $r$ on $(X,g,\xi)$, which is stable. Then there exists $\epsilon>0$ such that 
    \[
    \mu_{\Sas}(E_\xi)-\mu_{\Sas}(\mathcal F)\geq \epsilon
    \]
    for any saturated transverse coherent subsheaf $\mathcal F\subset E_\xi$ of $0<\rank \mathcal F<r$.
\end{Thm}
\begin{proof}
    Throughout this proof, all $L^2$ norms are taken with respect to the measure $d\eta^n\wedge \eta$. Let $s=\rank \mathcal F$. 
    By the Kobayashi-Hitchin correspondence, there is a $T$-invariant transverse Hermitian--Einstein metric $h$ on $E_\xi$. Take an arbitrary saturated transverse coherent subsheaf $\mathcal F\subset E_\xi$. Let $Z\subset X$ be the singular locus of $\mathcal F\subset E$. Then $Z$ is a transverse analytic subvariety of transverse complex codimension at least $2$. Let $\pi: E|_{X-Z}\to E|_{X-Z}$ be the $h$-orthogonal projection associated with subbundle $\mathcal F|_{X-Z}\subset E|_{X-Z}$. Since both $\mathcal F$ and $h$ are $T$-invariant, the orthogonal projection $\pi$ is $T$-invariant, and hence basic, i.e. 
    \[
    \pi\in A^0_B(X-Z,\End(E)),
    \]
    where $A^0_B(X-Z,\End(E))$ is the space of basic smooth sections of $\End(E)$. Let $h_F$ be the restriction of $h$ to $\mathcal F|_{X-Z}$. Then we have 
    \[
    \sqrt{-1}\Lambda\tr(F_{h_{F}})=\sqrt{-1}\Lambda \tr(\pi F_h\pi)-|\nabla^{1,0}_{\End(E),B}\pi|^2,
    \]
    where $\nabla_B=\nabla^{1,0}_B+\nabla^{0,1}_B$ is the transverse covariant derivative. Note that $|\nabla^{0,1}_{\End(E),B}\pi|=|\nabla^{1,0}_{\End(E),B}\pi|$ since $\pi$ is self-adjoint. Integrating this, we have 
    \[
    \mu_{\Sas}(E)-\mu_{\Sas}(F)=\frac{1}{4\pi ns \vol(X)}\int_{X-Z}|\nabla^{1,0}_B\pi|^2(d\eta)^n\wedge \eta,
    \]
    Thus $\nabla_B\pi\in L^2_B(X-Z,\End(E))$. Since $|\pi|^2=\tr(\pi)=s$, $\pi\in L^\infty_B(X-Z,\End(E))$. Combining this, we have 
    \[
    \pi\in W^{1,2}_B(X-Z,\End(E)).
    \]
    Since $Z$ is transverse analytic subvariety of transverse complex codimension at least $2$, its $W^{1,2}$-capacity is $0$. Indeed, taking a foliated chart $\phi: I\times B^{2n}\to X$, we have a cut-off function $\chi_\epsilon$, ($\epsilon>0$) obtained by smoothing 
    \[
    1-\max\{0,1-\epsilon^{-1}(d(\cdot ,Z)-\epsilon)_+\},
    \]
    which does not depend on the Reeb direction, such that 
    \[
    \int|d\chi_\epsilon|^2=O(\epsilon^2)\to 0.
    \]
    Since $|\pi|^2=s$ is constant and $\pi$ is self-adjoint, we obtain
    \[
    {\rm Re}\langle \pi, \nabla \pi\rangle_h=0.
    \]
    Hence, 
    \[
    \int_X|\nabla (\chi_\epsilon\pi)|^2\leq s\cdot \int_X|d\chi_\epsilon|^2 +\int_{X-Z}|\nabla\pi|^2.
    \]
    Furthermore, we have $\chi_\epsilon\pi\to \pi$ in $L^2$ by dominated convergence theorem. Thus there exists a subsequence $\epsilon_j\to 0$ such that $\chi_{\epsilon_j}\pi$ converges weakly in $W^{1,2}_B(X,\End(E))$. Since its weak $W^{1,2}$-limit must agree with $L^2$-limit, we conclude that
    \[
    \pi\in W^{1,2}_B(X,\End(E)).
    \]
    
    Thus, it suffices to show that there is uniform constant $\lambda>0$ such that 
    \begin{equation}\label{spectral gap}
        \|\nabla^{0,1}_Bu\|_{L^2}^2\geq \lambda \|u\|_{L^2}^2
    \end{equation}
    for any self adjoint $u\in W^{1,2}_B(\End(E))$ with $\tr(u)=0$. Indeed, if this inequality holds, we have 
    \[
    \|\pi-\frac{s}{r}I_E\|_{L^2}^2=\frac{s(r-s)}{r}\vol(X)
    \]
    and thus 
    \[
    \|\nabla^{0,1}_B\pi\|_{L^2}^2\geq \lambda\cdot \frac{s(r-s)}r\vol(X),
    \]
    and hence 
    \[
    \mu_{\Sas}(E)-\mu_{\Sas}(\mathcal F)\geq \frac{r-s}{4\pi nr}\lambda.
    \]
    Thus taking $\epsilon:=\frac{\lambda}{4\pi nr}$, we get the required property. If there is no such $\lambda$ satisfying the inequality (\ref{spectral gap}), we have a sequence $(u_j)\in W^{1,2}_B(\End(E))$ such that all $u_j$ are self-adjoint and 
    \[
    \|u_j\|_{L^2}=1,\  \|\nabla^{0,1}_Bu_j\|_{L^2}\to 0,\ \tr(u_j)=0.
    \]
    Since $\|\nabla^{1,0}_Bu_j\|_{L^2}=\|\nabla^{0,1}_Bu_j\|$ by self-adjointness, we have
    \[
    \|u_j\|_{1,2}^2=\|u_j\|_{L^2}^2+2\|\nabla^{0,1}_Bu_j\|_{L^2}^2,
    \]
    and hence $(u_j)$ is a bounded sequence in $W^{1,2}_B(\End(E))$. Now, there is an embedding $W^{1,2}_B(X,\End(E))\to W^{1,2}(X,\End(E))$, since $\nabla=\nabla_B+\eta\otimes \mathcal L_\xi$ and hence $\nabla u=\nabla_Bu$ for any basic smooth section. Moreover, since basic smooth functions are $T$-invariant by density of the Reeb orbit and $T$-action is isometric, we have $W^{1,2}_B=(W^{1,2})^T$, especially, it is a closed subspace. Thus, by the Rellich compactness theorem and Banach--Alaoglu theorem, after passing to a subsequence, we may assume that there is $u\in W^{1,2}_B(\End(E))$ such that $u_j$ converges to $u$ in the strong topology in $L^2$ and the weak topology in $W^{1,2}_B$. Then we have 
    \[
    \|u\|_{L^2}=1, \ u=u^*,\ \tr(u)=0,\ \nabla^{0,1}_Bu=0.
    \]
    Note that the last equality is at first in the sense of weak solution. However, this equation shows that $(\nabla^{0,1}_B)^*\nabla^{0,1}_Bu=0$ and the same argument as in the proof of the basic elliptic regularity statement \cite[Lemma 2.8]{BARAGLIA2022108661} shows that $u$ is smooth and $\nabla^{0,1}u=0$ in the strong sense. Now, $\mathcal L_\xi u$ is also $0$, and hence $u\in H^0_B(X,\End^0(E))$. However, $E_\xi$ is simple by \cite[Proposition 3.24]{BARAGLIA2022108661} and thus $u=0$. This contradicts $\|u\|_{L^2}=1$.
\end{proof}
We can now prove Theorem \ref{Openness of Stability I}.

\noindent
\begin{proof}[Proof of Theorem \ref{Openness of Stability I}.]
    We divide the proof into two steps.
    
    \noindent
    \underline{Step 1}. In this step, we prove that there is a neighbourhood $U$ of $\tau_0$ such that for any saturated $\mathbb X^*(T_\mathbb C)$-graded module $N\subset M$, $\mu_{\gr}(N_\tau)<\mu_{\gr}(M_\tau)$. By Theorem \ref{Spectral Gap Theorem} and Theorem \ref{equivalence of stabilities}, there exists $\epsilon>0$, which doesn't depend on $N$, such that 
    \[
    \mu_{\gr}(N_{\tau_0})\leq \mu_{\gr}(M_{\tau_0})-\epsilon.
    \]
    By Lemma \ref{First derivative estimate}, we have
    \[
    \|d\mu(N_\tau)\|\leq C_M+C'_\tau(C_M\|\tau\|-\mu(N_\tau)).
    \]
    Assume that there is a sequence $\tau_i\to \tau_0$ with saturated graded modules $N^{(i)}\subset M$ such that 
    \[
    \mu(N^{(i)}_{\tau_i})\geq \mu(M_{\tau_i}).
    \]
    Since $M_{\tau_0}$ is stable, we have 
    \[
    \mu(N^{(i)}_{\tau_0})\leq \mu(M_{\tau_0})-\epsilon.
    \]
    Let $\gamma_i(t)=(1-t)\tau_0+t\tau_i$ and set 
    \[
    f_i(t)=\mu(N^{(i)}_{\gamma_i(t)})-\mu(M_{\gamma_i(t)}).
    \]
    Then we have $f_i(0)\leq -\epsilon$ and $f_i(1)\geq 0$. Since $f_i$ is smooth, we can choose $t_i,t'_i$ ($t_i<t'_i$) such that $f_i(t_i)=-\frac{2}{3}\epsilon$, $f_i(t'_i)=-\frac {1}{3}\epsilon$, and 
    \begin{equation}\label{how to take $t_i$ and $t'_i$}
        f_i(t)\in \left[-\frac{2}{3}\epsilon,-\frac{1}{3}\epsilon\right],
    \end{equation} 
    for any $t\in [t_i,t'_i]$. By the mean value theorem, we have some $u_i\in [t_i,t'_i]$ such that
    \[
    f_i'(u_i)=\frac{\epsilon}{3(t'_i-t_i)}\geq \frac{\epsilon}{3}.
    \]
    Let $\zeta_i=\gamma_i(u_i)$. Then we have
    \[
    \|d\mu(N^{(i)}_{\zeta_i})-d\mu(M_{\zeta_i})\|\geq \frac{\epsilon}{3\|\tau_0-\tau_i\|}\to \infty  \qquad (i\to \infty).
    \]
    Since $\mu(M_{\tau})$ is smooth in $\tau$, we have 
    \[
    C_M+C'_{\zeta_i}(C_M\|\zeta_i\|-\mu(N^{(i)}_{\zeta_i}))\geq \|d\mu(N^{(i)}_{\zeta_i})\|\to \infty.
    \]
    Since $C'_{\zeta_i}\to C'_{\tau_0}$ and $\zeta_i\to \tau_0$, we have 
    \[
    \mu(N^{(i)}_{\zeta_i})\to -\infty \qquad (i\to \infty).
    \]
    This contradicts the equation (\ref{how to take $t_i$ and $t'_i$}), since $f_i(u_i)$ and $\mu(M_{\zeta_i})$ are bounded.

    \noindent
    \underline{Step 2}. Take an open neighbourhood $U\subset \mathfrak t_+$ of $\tau_0$ obtained in Step 1. Fix $\tau \in U$. If the Reeb flow of $\tau$ is dense in $T$, $E_\tau$ is stable by Theorem \ref{subsheaf-submodule correspondence} and Theorem \ref{equality of slopes}. Thus we may assume that the Reeb flow generated by $\tau$ is not dense in $T$. Take an arbitrary saturated submodule $N_\tau\subset M_\tau$. Since $M_\tau\simeq M\simeq M_\xi$ as ungraded $R$-module, we can endow the filtration $(F^\lambda M)_{\lambda\in \mathbb R}$of $M$ such that 
    \[
    F^\lambda M=\bigoplus_{\langle\chi,\xi\rangle\geq\lambda}M_\chi,
    \]
    This filtration is induced by the valuative function $v_\xi$, and its associated graded module is isomorphic to $M_\xi$. Let $(F^\lambda N_\tau)$ be an induced filtration and let $N':=\gr^{F}(N_\tau)\subset \gr^F(M)$. Then $N'$ is a graded submodule of $M_\xi$, and hence it can be regarded as a graded submodule of $M$. For each graded piece $N_{\tau,\lambda}$, this filtration induces a bounded filtration of a finite-dimensional vector space. Thus, taking graded object, we have 
    \[
    \dim(N_{\tau,\lambda})=\sum_{\langle\chi,\tau\rangle=\lambda}\dim(N'_\chi)=\dim(N'_{\tau,\lambda}).
    \]
    Thus $N_\tau$ and $N'_\tau$ have the same Hilbert function, especially, we have
    \[
    \mu(N_\tau)=\mu(N'_{\tau}),\  \rank(N_\tau)=\rank(N'_\tau)
    \]
    Since taking saturation does not decrease the slope, we have 
    \[
    \mu(N_\tau)\leq \mu((N')^{{\rm sat}}_\tau)
    \]
    Since $\tau\in U$, $\mu((N')^{{\rm sat}}_\tau)<\mu(M_\tau)$ by Step 1. Therefore, $M_\tau$, and hence $E_\tau$ is stable.
\end{proof}

In the argument of the Step 2 of the proof, we have the following important corollary:
\begin{algcorollary}
    Let $E$ be a $T$-linearized transverse holomorphic vector bundle on a compact Sasakian manifold $(X,g,\xi)$ of dimension $2n+1$, where $n\geq 1$. Then for any $\tau\in \mathfrak t_+$, any saturated subsheaf $\mathcal F_\tau\subset E_\tau$ admits a saturated $T$-invariant subsheaf $\mathcal F'_\tau\subset E_\tau$ with the same rank such that 
    \[
    \mu(\mathcal F_\tau)\leq \mu(\mathcal F'_\tau).
    \]
    In particular, the maximal destabilising subsheaf $\mathcal F_\tau\subset E_\tau$ is $T$-invariant.
\end{algcorollary}
\begin{proof}
    Let $N_\tau\subset M_\tau$ be a saturated graded submodule. We can construct a saturated $\mathbb X^*(T_\mathbb C)$-graded submodule $N'\subset M$ with $\mu(N_\tau)\leq \mu(N'_\tau)$, and $\rank(N_\tau)=\rank(N'_\tau)$, by the construction in Step 2 in the proof of Theorem \ref{Openness of Stability I}. If $N_\tau$ is maximal destabilising, we have $N_\tau=N'_\tau$ by the uniqueness of the maximal destabilising subsheaf \cite[Theorem 2.4]{hada2026optimalalgebraictangentcone}. Together with Theorem \ref{subsheaf-submodule correspondence}, this proves the desired statement.
\end{proof}

\begin{Rem}
    This may be viewed as a transverse analogue of the standard fact that the Harder--Narasimhan filtration of an equivariant sheaf is equivariant; see \cite[Proposition 4.13]{KOOL20111700}, and \cite[Theorem 2.1]{BDGP21} for the case of projective toric varieties.
\end{Rem}

Since $\mu(E_{\lambda\tau})=\lambda\mu(E_\tau)$ for any $\tau\in \mathfrak t_+$ and $\lambda>0$, the set 
\[
\mathfrak t_+^{\st}(E)=\mathfrak t_+^{\st}:=\{\tau\in \mathfrak t_+\mid E_\tau \text{ is stable}\}
\]
is an open cone. We call it the {\bf stable cone} of $E$. When the cone $\mathcal C=\Spec(R)$ is $\mathbb Q$-factorial, the consequence is simpler. Indeed, we have the following lemma:
\begin{alglemma}\label{liniarity of slopes when cone is Q-factorial}
    Let $\mathcal C$ be a normal $\mathbb Q$-factorial affine variety with good $T_\mathbb C$-action. 
    Then, for any reflexive $\mathbb X^*(T_\mathbb C)$-graded $R$-module $M$,
    \[
    \mu(M): \mathfrak t_+\to \mathbb R;\ \tau\mapsto \mu(M_\tau)
    \]
    is linear. As a consequence, the stable cone $\mathfrak t_+^{\st}(E)$ is convex for any $T$-linearized transverse holomorphic vector bundle.
\end{alglemma}
\begin{proof}
    By Lemma \ref{determinant degree equals the initial degree}, we have $\rank (M)\cdot \mu(M_\tau)=\mu(\det(M)_\tau)$ for any $\tau\in \mathfrak t_+$. Let $r>0$ be a positive integer such that $L:=\det(M)^{[r]}$ is free at the origin. 
    Since $L_{\mathfrak m_0}$ is a free $R_{\mathfrak m_0}$-module of rank $1$, $\dim(L/\mathfrak m_0L)=1$. Let $s\in L$ be a homogeneous element which defines a generator of $L/\mathfrak m_0L$. Then we have a morphism
    \[
    s: R(-\deg(s))\to L,
    \]
    which is an isomorphism at the vertex $0$. Hence the support of its cokernel is $T_\mathbb C$-invariant closed set, which does not contain the vertex $0$. Since $T_\mathbb C$-action is good, we conclude that the above morphism is an isomorphism. Thus we have 
    \[
    \det(M)^{[r]}=R(\chi)
    \]
    for some $\chi\in \mathbb X^*(T_\mathbb C)$. Since $\det(M)^{[r]}\simeq \det(M^{\oplus r})$, we have 
    \[
    \deg(M_\tau)=\frac{1}{r}\deg((M^{\oplus r})_\tau)=\frac{1}{r}\deg(\det(M)^{[r]}_\tau)=\frac{a_{n+1}(R_\tau)}{r}\langle\chi,\tau\rangle.
    \]
    Thus we have 
    \[
    \mu(M_\tau)=\frac{1}{r\cdot \rank M}\langle\chi,\tau\rangle.
    \]
    This is obviously linear. 
    
    Let $E$ be a $T$-linearized transverse holomorphic vector bundle, and $M$ an associated $\mathbb X^*(T_\mathbb C)$-graded $R$-module. Then, $\mathfrak t_+^{\st}:=\mathfrak t_+^{\st}(E)$ can be written as follows:
    \[
    \mathfrak t_+^{\st}=\bigcap_{N\subset M}\{\tau\in \mathfrak t_+\mid \mu(M_\tau)-\mu(N_\tau)>0\},
    \]
    where the intersection runs over all saturated $\mathbb X^*(T_\mathbb C)$-graded submodules. By the linearity of $\mu$, we have that $W_N:=\{\tau\in \mathfrak t_+\mid \mu(M_\tau)-\mu(N_\tau)>0\}$ is the intersection of an open half-space with $\mathfrak t_+$. By construction, $\mathfrak t_+$ is convex and hence $\mathfrak t_+^{\st}$ is also convex.
\end{proof}
We also have a notion of {\bf semistable cone}, i.e. 
\[
\mathfrak t_+^{{\rm ss}}=\{\tau\in \mathfrak t_+\mid E_\tau\ \text{is semistable}\}=\bigcap_{N\subset M}\{\tau\in \mathfrak t_+\mid \mu(M_\tau)-\mu(N_\tau)\geq 0\},
\]
where $M$ is the corresponding $\mathbb X^*(T_\mathbb C)$-graded module and $N$ runs over all saturated graded submodules of $0<\rank(N)<\rank(M)$. By definition, this is a closed cone in $\mathfrak t_+$. Obviously, we have 
\[
\overline{\mathfrak t_+^{\st}}^{\mathfrak t_+}\subset \mathfrak t_+^{{\rm ss}}
\]
\begin{prop}
    Let $(X,g,\xi)$ be a compact Sasakian manifold of $\dim(X)=2n+1$, $n\geq 1$, whose affine cone is $\mathbb Q$-factorial. Let $E_\xi$ be a transverse holomorphic vector bundle. Assume that $\mathfrak t_+^{\st}(E)$ is non-empty. Then we have 
    \[
    \overline{\mathfrak t_+^{\st}(E)}^{\mathfrak t_+}=\mathfrak{t}^{{\rm ss}}_+.
    \]
\end{prop}
\begin{proof}
    Fix $\tau_*\in \mathfrak t_+^{\st}(E)$. Take an arbitrary $\tau\in \mathfrak t_+^{{\rm ss}}$. Then we have 
    \[
    \mu(M_\tau)-\mu(N_\tau)\geq 0, \mu(M_{\tau_*})-\mu(N_{\tau_*})>0
    \]
    for any saturated $\mathbb X^*(T_\mathbb C)$-graded submodules $N$. Since $\mu(M_\tau)$ and $\mu(N_\tau)$ is linear on $\mathfrak t_+$, we have 
    \[
    \mu(M_{(1-t)\tau+t\tau_*})-\mu(N_{(1-t)\tau+t\tau_*})>0
    \]
    for any $N$ and $0<t\leq 1$. Thus $\tau\in \overline{\mathfrak t_+^{\st}(E)}^{\mathfrak t_+}$. 
\end{proof}
\subsection{Openness and finiteness: an algebro-geometric proof}
In this section, we mainly treat the following algebraic case:
\begin{itemize}
    \item $\mathcal C$ is a normal affine variety of dimension $n+1$ for $n\geq 1$ with good $T_\mathbb C$-action of rank $k$.
    \item $\mathcal E$ is a $T_\mathbb C$-linearized torsion-free sheaf.
\end{itemize}
Let $R:=\Gamma(\mathcal C,\mathcal O_{\mathcal C})$ be the ring of regular functions. The $T_\mathbb C$-action on $\mathcal C$ induces the weight decomposition
\[
R=\bigoplus_{\vec m\in \mathbb X^*(T_\mathbb C)}R_{\vec m}
\]
and $M:=\Gamma(\mathcal C,\mathcal E)$ is a $\mathbb X^*(T_\mathbb C)$-graded $R$-module. Let 
\[
\omega:=\langle\vec m\in \mathbb X^*(T_\mathbb C)\mid R_{\vec m}\neq 0\rangle_{\mathbb R_{\geq 0}}
\]
be a weight cone and $\sigma\subset \mathbb X_*(T_\mathbb C)_{\mathbb R}$ be the dual of $\omega$.

By \cite[Theorem 3.4]{AH06}, after resolving the singularities if necessary, we have a smooth semiprojective variety $Y$ and $\sigma$-pp-divisor $\mathfrak D$ such that
\[
R_{\vec m}=H^0(Y, \mathcal O_Y(\lfloor \mathfrak D(\vec m)\rfloor)),\qquad (\forall \vec m\in \omega\cap \mathbb X^*(T_\mathbb C)).
\]
Since $T_\mathbb C$-action on $\mathcal C$ is good, $R_0=H^0(Y,\mathcal O_Y)=\mathbb C$. Therefore, $Y$ is indeed projective.
$T_\mathbb C$-invariant prime divisors on $\mathcal C$ are of two types (See \cite[Remark 3.16]{PS11}):
\begin{enumerate}
    \item There is some prime divisor $Z$ on $Y$ and a vertex $v\in \Delta_Z$, which is the $Z$-coefficient of $\mathfrak D$, such that corresponding ideal is defined by 
    \[
    (I_{Z,v})_{\vec m}=R_{\vec m}\cap \{f\in \mathbb C(Y)\mid v_Z(f)>-\langle \vec m,v\rangle\}
    \]
    Furthermore, in this case, $\mathfrak D(u)|_Z$ is big for any $u\in {\rm Int}((\Delta_Z-v)^\vee)$.
    \item There is a ray $\rho$ of $\sigma$ such that corresponding ideal is 
    \[
    (I_\rho)_{\vec m}=\begin{cases}
        R_{\vec m}&(\vec m\notin \rho^\perp)\\
        0&(\vec m\in \rho^\perp).
    \end{cases}
    \]
    Furthermore, in this case, $\mathfrak D(u)$ is big for some (equivalently, for any) $u\in {\rm reint }(\omega\cap \rho^\perp)$.
\end{enumerate}

Note that in the first case, all but finitely many $Z$ satisfy $\Delta_Z=\sigma$. In this case, $v$ have to be $0$, and we have 
\[
(I_{Z,v})_{\vec m}=H^0(Y,\mathcal O_Y(\lfloor \mathfrak D(\vec m)\rfloor -Z))
\]
and thus we have the following lemma:
\begin{alglemma}\label{restricted volume lemma}
    Let $Z\subset Y$ be a prime divisor with $\Delta_Z=\sigma$ and $\mathfrak D(u)|_Z$ is big for every $u\in {\rm Int}(\sigma^\vee)$. Then we have 
    \[
    a_n((R/I_{Z,0})_\tau)=\frac{1}{(\dim(Y)-1)!}\int_{\sigma^\vee}e^{-\langle x,\tau\rangle}\vol_{Y|Z}(\mathfrak D(x))dx,
    \]
    where $dx$ is the lattice normalized measure.
\end{alglemma}
\begin{proof}
    Put $d=\dim(Y)$ and $k=\dim(T_\mathbb C)$, so that $n=d+k-1$, and $Q=R/I_{Z,0}$. Take a finite rational simplicial subdivision of $\sigma^\vee$ such that $\mathfrak D$ is linear on each cone. Since the contribution of proper faces is of lower order, it is enough to work on a fixed full-dimensional rational simplicial cone. From now on, assume that $\sigma^\vee$ is already simplicial and $\mathfrak D$ is linear on $\sigma^\vee$.

    Let $\vec m_1,\dots,\vec m_k$ be the primitive ray generators of $\sigma^\vee$, and put $D_i:=\mathfrak D(\vec m_i)$. Choose $e>0$ such that all $eD_i$ are base-point-free Cartier divisors. They define a morphism
    \[
    \phi: Y\to \Pi:=\mathbb P(H^0(Y,eD_1))\times \dots\times \mathbb P(H^0(Y,eD_k)).
    \]
    Let $W$ be the image of $\phi$, and let $L_i$ be restrictions of $p_i^*\mathcal O_{\mathbb P(H^0(Y,eD_i))}(1)$, where $p_i: \Pi\to \mathbb P(H^0(Y,eD_i))$ are $i$-th projection morphisms. Then we have 
    \[
    \phi^*\mathcal O_W(a_1L_1+\dots+a_kL_k)=\mathcal O_Y(\sum ea_iD_i).
    \]
    Since every positive combination of the $D_i$'s is big, $\dim(W)=d$. Moreover, by the assumption that $\mathfrak D(u)|_Z$ is big for $u\in {\rm Int}(\sigma^\vee)$, we have $\dim(\phi(Z))=d-1$. 

    Set 
    \[
    \Lambda_+=\sum_{i=1}^ke\mathbb Z_{\geq 0}\cdot \vec m_i\subset \sigma^\vee \cap \mathbb X^*(T_\mathbb C).
    \]
    Then $\Lambda_+\subset \sigma^\vee \cap \mathbb X^*(T_\mathbb C)$ is sub-semigroup such that $\Lambda=\langle\Lambda_+\rangle_\mathbb Z\subset \mathbb X^*(T_\mathbb C)$ is the sublattice with finite index. For any $\vec m\in \sigma^\vee \cap \mathbb X^*(T_\mathbb C)$, we can write
    \[
    \vec m=p(\vec m)+\delta(\vec m),
    \]
    with $p(\vec m)\in \Lambda_+$ and $\delta(\vec m)\in (\sum_i [0,e) \vec m_i)\cap \mathbb X^*(T_\mathbb C)$, which is a finite set. Put
    \[
    D_\delta=\lfloor \mathfrak D(\delta)\rfloor, \qquad \mathcal G_\delta:={\rm Cok}(\phi_*(\mathcal O_Y(D_\delta-Z))\to \phi_*(\mathcal O_Y(D_\delta)))\subset \phi_*\mathcal (O_Y(D_\delta)|_Z).
    \] 
    Since $\mathfrak D(p(\vec m))$ is integral Cartier, we have 
    \[
    \lfloor \mathfrak D(\vec m)\rfloor=\mathfrak D(p(\vec m))+D_{\delta(\vec m)}.
    \]
    By Fujita's vanishing theorem, we have $m_0>0$, independent of $\delta$, such that 
    \[
    H^p\left(W,\mathcal O_W(\sum_{i=1}^k a_iL_i)\otimes \phi_*\mathcal O_Y(D_\delta-Z)\right)=H^p\left(W,\mathcal O_W(\sum_{i=1}^k a_iL_i)\otimes \phi_*\mathcal O_Y(D_\delta)\right)=0
    \]
    for all $p>0$, $(a_i)\subset \mathbb Z$ with $\min\{a_i\}\geq m_0$. Thus 
    \[
    \dim(Q_{\vec m})=h^0(W,\mathcal O_W(a_iL_i)\otimes \mathcal G_\delta),
    \]
    for $\vec m\in \sigma^\vee \cap \mathbb X^*(T_\mathbb C)$ satisfying $p(\vec m)=\sum_i ea_i\vec m_i$ with $a_i\geq m_0$. 
    The right hand side is given by the multigraded Hilbert polynomial. Since $\phi|_Z$ is generically finite, the top-dimensional cycle of $\mathcal G_\delta$ is $\phi_*[Z]$. Indeed, we have $\dim(W)=\dim(Y)$ and $\phi|_Z$ is generically finite. Hence the pushforward of
    \[
    0\to \mathcal O_Y(D_\delta-Z)\to \mathcal O_Y(D_\delta)\to \mathcal O_Z(D_\delta|_Z)\to 0
    \] 
    remains exact near the generic point of $\phi(Z)$. Therefore, we have 
    \[
    \mathcal G_\delta\simeq \phi_*\mathcal O_Z(D_\delta|_Z)
    \]
    near the generic point of $\phi(Z)$, and consequently, 
    \[
    [\mathcal G_\delta]=\phi_*[Z].
    \]
    Hence its homogeneous part of degree $d-1$ is
    \[
    \frac{1}{(d-1)!}\left(\sum_i a_iL_i\right)^{d-1}\cdot \phi_*[Z]=\frac{1}{(d-1)!}\mathfrak D(p)^{d-1}\cdot Z=\frac{1}{(d-1)!}\vol_{Y|Z}(\mathfrak D(p)).
    \]
    Since the set of possible $\delta$ is finite, we have some constant $C>0$ such that 
    \[
    \left|\dim(Q_{\vec m})-\frac{\vol_{Y|Z}(\mathfrak D(p(\vec m)))}{(d-1)!}\right|\leq C(1+\|\vec m\|)^{d-2},
    \]
    whenever $a_i\geq m_0$. For each fixed $\delta$,
    \[
    \mathfrak D(\delta+p)^{d-1}\cdot Z-\mathfrak D(p)^{d-1}\cdot Z
    \]
    is polynomial in $p$ of degree at most $d-2$ if $d\geq 2$ and $0$ if $d=1$. Since the possible $\delta$ is finite, we can re-choose the constant $C>0$ so that 
    \begin{equation}\label{strong estimate in enough interiour of sigma}
    \left|\dim(Q_{\vec m})-\frac{\vol_{Y|Z}(\mathfrak D(\vec m))}{(d-1)!}\right|\leq C(1+\|\vec m\|)^{d-2},
    \end{equation}
    whenever $a_i\geq m_0$. In particular, we have 
    \[
    \left|t^{d-1}\dim(Q_{\vec m})-\frac{\vol_{Y|Z}(\mathfrak D(t\vec m))}{(d-1)!}\right|\leq Ct(t+\|t\vec m\|)^{d-2},
    \]
    whenever $a_i\geq m_0$. Here we used the $(d-1)$-homogeneity of $\vol_{Y|Z}$. Let $\Gamma_1$ be the set of such lattice points, and $\Gamma_2$ be its complement in $\sigma^\vee \cap \mathbb X^*(T_\mathbb C)$. $\Gamma_2$ is contained in a finite union of bounded translates of codimension one faces.
    Moreover, we have
    \[
    \dim Q_{\vec m}\leq h^0(Z,\mathcal O_Z(\lfloor \mathfrak D(\vec m)\rfloor))\leq C_2(1+\|\vec m\|)^{\dim(Y)-1}
    \]
    for some constant $C_2>0$ for all $\vec m\in \sigma^\vee \cap \mathbb X^*(T_\mathbb C)$. It follows that
    \[
    t^n\sum_{\vec m\in \Gamma_2}e^{-\langle t\vec m,\tau\rangle}\dim(Q_{\vec m})=O(t),\qquad (t\to 0).
    \]
    Using (\ref{strong estimate in enough interiour of sigma}), we have 
    \begin{align*}
    t^nF_Q(\tau,t):=&t^n\sum_{\vec m\in \sigma^\vee \cap \mathbb X^*(T_\mathbb C)}e^{-\langle t\vec m,\tau\rangle}\dim(Q_{\vec m})\\
    =&t^n\sum_{\vec m\in \Gamma_1}e^{-\langle t\vec m,\tau\rangle}\dim(Q_{\vec m})+O(t)\\
    =&\frac{t^k}{(d-1)!}\sum_{\vec m\in \Gamma_1}e^{-\langle t\vec m,\tau\rangle}\vol_{Y|Z}(\mathfrak D(t\vec m))+O(t).
    \end{align*}
    The sum on the last line is a Riemann sum, and hence converges to 
    \[
    \frac{1}{(d-1)!}\int_{\sigma^\vee}e^{-\langle x,\tau\rangle}\vol_{Y|Z}(\mathfrak D(x))dx.
    \]
\end{proof}
The following corollary follows from this lemma.
\begin{algcorollary}\label{finiteness corollary}
    Let $\mathcal P_1$ be the set of $T_\mathbb C$-invariant prime divisors $P=D_{Z,0}$ of type (i) on $\mathcal C$ with $\Delta_Z=\sigma$. Then there exist finitely many positive functions
    \[
    \Phi_1,\dots,\Phi_N: \mathfrak t_+={\rm Int}(\sigma)\to \mathbb R
    \]
    and a map
    \[
    q: \mathcal P_1\to \mathbb Z_{\geq 0}^N
    \]
    such that 
    \[
    \Theta_P(\tau)=\sum_{j=1}^N q_j(P)\Phi_j(\tau).
    \]
\end{algcorollary}
\begin{proof}
    Choose a finite rational simplicial subdivision $\Sigma$ of $\sigma^\vee$ such that $\mathfrak D$ is linear for each cone $\kappa\in \Sigma$. For $\rho\in \Sigma(1)$, let $\vec m_\rho$ be its primitive generator. Since $\Sigma(1)$ is finite and $\mathfrak D(u)$ is semiample for any $u$, we can take $e\in \mathbb Z_{>0}$ so that 
    \[
    A_{\rho}:=e\mathfrak D(\vec m_\rho)
    \]
    is base-point-free Cartier divisor for any $\rho\in \Sigma(1)$. Let $\kappa\in \Sigma$ be a maximal cone of $\Sigma$, with rays $\rho_{1,\kappa},\dots,\rho_{k,\kappa}$. For $x\in \kappa$, we have the unique $a_{i,\kappa}(x)\in \mathbb R_{\geq 0}$ such that 
    \[
    x=\sum_{i=1}^k a_{i,\kappa}(x)\vec m_{i,\kappa},
    \]
    where $\vec m_{i,\kappa}=\vec m_{\rho_{i,\kappa}}$. Then we have 
    \[
    \mathfrak D(x)=\sum_{i=1}^k \frac{1}{e}a_{i,\kappa}(x)A_{i,\kappa}.
    \]
    Since $\mathfrak D(x)$ is semiample, we have 
    \[
    \vol_{Y|Z}(\mathfrak D(x))=\mathfrak D(x)^{d-1}\cdot Z=\frac{1}{e^{d-1}}\sum_{|\alpha|=d-1}\binom{d-1}{\alpha}a_\kappa^\alpha q_{\kappa,\alpha}(Z),
    \]
    where 
    \[
    q_{\kappa,\alpha}(Z)=(A_{1,\kappa}^{\alpha_1}\cdot \dots\cdot A_{k,\kappa}^{\alpha_k}\cdot Z)\in \mathbb Z_{\geq 0}.
    \]
    Thus, by Lemma \ref{restricted volume lemma}, we have 
    \[
    \Theta_P(\tau)=\sum_{\kappa,\alpha}\frac{q_{\kappa,\alpha}(Z)}{e^{d-1}\alpha!a_{n+1}(R_\tau)}\int_{\kappa}e^{-\langle x,\tau\rangle}a_\kappa(x)^\alpha dx.
    \]
    Let 
    \[
    \Phi_{\kappa,\alpha}(\tau)=\frac{1}{e^{d-1}\alpha!a_{n+1}(R_\tau)}\int_{\kappa}e^{-\langle x,\tau\rangle}a_{\kappa}(x)^\alpha dx.
    \]
    Then $\Phi_{\kappa,\alpha}(\tau)>0$ and we have 
    \[
    \Theta_P(\tau)=\sum_{\kappa,\alpha}q_{\kappa,\alpha}(Z)\cdot \Phi_{\kappa,\alpha}(\tau).
    \]
\end{proof}
Let $\mathcal P_2$ be the set of all $T_\mathbb C$-invariant prime divisors on $\mathcal C$ which is {\it not} contained in $\mathcal P_1$. Since $\Delta_Z=\sigma$ for all but finite $Z$ and divisors of type (ii) are finite, $\mathcal P_2$ is a finite set. 

\begin{algtheorem}\label{finiteness of maximal slope functions}
    For any $T_\mathbb C$-graded torsion-free module $M$, let $\mathfrak F=\mathfrak F_M$ be the family of all saturated $\mathbb X^*(T_\mathbb C)$-graded submodules, and $F_s$ ($s=1,\dots,\rank (M)-1$) be the set defined by 
    \[
    F_s=\{\mu(N): \mathfrak t_+\to \mathbb R\mid N\in \mathfrak F, \rank (N)=s\}.
    \]
    Then, for any $s\in \{1,\dots,\rank(M)-1\}$, the set $F^{\max}_s$ of functions in $F_s$ defined by the following property is a finite set:
    \begin{itemize}
        \item $f\in F_s^{\max}$ if and only if $f\in F_s$ and there is no $g\neq f\in F_s$ such that $g\geq f$ on $\mathfrak t_+$.
    \end{itemize}
\end{algtheorem}
\begin{proof}
    For each $s\in \{1,\dots,\rank(M)-1\}$, embed
    \[
    \left(\bigwedge^sM\right)^{**}\hookrightarrow \bigoplus_{i=1}^{k_s}R(\beta_{s,i})
    \]
    for some $\beta_{s,i}\in \mathbb X^*(T_\mathbb C)$. Let $\mathfrak F_s$ be a family of all saturated $\mathbb X^*(T_\mathbb C)$-graded submodules of rank $s$ of $M$. For $f,g\in F_s$, we say $f\leq g$ if $f(\tau)\leq g(\tau)$ for any $\tau\in \mathfrak t_+$. We prove that $F_s$ has only finitely many maximal elements under the assumption of the statement. 
    
    For any saturated $\mathbb X^*(T_\mathbb C)$-graded submodule $N\subset M$ of rank $s$, we have an inclusion
    \[
    \det(N)\hookrightarrow R(\beta_{s,i})
    \]
    for some $i\in \{1,\dots,k_s\}$. Let $C$ be the cokernel of this map. Then we have 
    \[
    \mu(N)=\frac{1}{s}\beta_{s,i}-\frac{1}{s}\sum_{P}\ell_P(C)\Theta_P.
    \] 
    Let $\mathfrak F_{s,i}\subset \mathfrak F_s$ be the family of saturated submodules $N\subset M$ with $\rank(N)=s$ for which there is an inclusion $\det(N)\hookrightarrow R(\beta_{s,i})$ induced by $\det(N)\subset (\bigwedge^sM)^{**}$. Then the above argument shows that for any $N\in \mathfrak F_{s,i}$, we have 
    \[
    \mu(N)\leq \frac{1}{s}\beta_{s,i}.
    \]
    By Corollary \ref{finiteness corollary}, we have
    \begin{align*}
    \sum_{P}\ell_P(C)\Theta_P=&\sum_{P\in \mathcal P_1}\ell_P(C)\Theta_P+\sum_{P\in \mathcal P_2}\ell_P(C)\Theta_P\\
    =&\sum_{j=1}^N \left(\sum_{P\in \mathcal P_1}\ell_P(C)q_j(P)\right)\Phi_j+\sum_{P\in \mathcal P_2}\ell_P(C)\Theta_P.
    \end{align*}
    Define 
    \[
    q_j(C):=\sum_{P\in \mathcal P_1}\ell_P(C)q_j(P).
    \]
    Let $L=N+\#\mathcal P_2$ and define the function $\nu$ as
    \[
    \nu: \mathfrak F_{s,i}\to \mathbb Z_{\geq 0}^L; \nu(N)=((q_j(C))_{j=1}^N,(\ell_P(C))_{P\in \mathcal P_2}).
    \]
    For $v,w\in \mathbb Z_{\geq 0}^L$, we define $v\leq w$ iff $v_i\leq w_i$ for any $i$. For $N,N'\in \mathfrak F_{s,i}$, $\mu(N)\leq \mu(N')$ if $\nu(N)\geq \nu(N')$ and $\mu(N)=\mu(N')$ if $\nu(N)=\nu(N')$. Thus, for any $N\in \mathfrak F_{s,i}$ that defines maximal element in $F_{s,i}$, $\nu(N)$ is minimal in $\mathcal V_{s,i}:=\nu(\mathfrak F_{s,i})$ and the finiteness of minimal elements in $\mathcal V_{s,i}$ implies the finiteness of maximal elements in $F_{s,i}$. By Dickson's lemma, there are only finitely many minimal elements in $\mathcal V_{s,i}$. Thus the set
    \[
    F_{s,i}^{\max}:=\{\mu(N)\mid N\in \mathfrak F_{s,i}: \text{ maximal}\}
    \]
    is a finite set. Since $F_s^{\max}=\{\mu(N)\mid N\in \mathfrak F_s: \text{ maximal}\}$ is contained in the union of $F_{s,i}^{\max}$, $F_s^{\max}$ is also finite.
\end{proof}
This theorem yields the following theorem:
\begin{algtheorem}\label{finiteness of walls}
    Let $M$ be a torsion-free $\mathbb X^*(T_\mathbb C)$-graded $R$-module. Define the stable cone $\mathfrak t_+^{\st}(M)$ and semistable cone $\mathfrak t_+^{{\rm ss}}(M)$ as follows:
    \[
    \mathfrak t_+^{\st}(M):=\{\tau\in \mathfrak t_+\mid \text{$M_\tau$ is stable}\},\ \mathfrak t_+^{{\rm ss}}(M):=\{\tau\in \mathfrak t_+\mid \text{$M_\tau$ is semistable}\}.
    \]
    Then, the following properties hold:
    \begin{enumerate}
        \item $\mathfrak t_+^{\st}(M)$ is an open cone in $\mathfrak t_+$.
        \item There are finitely many saturated graded submodules $N_1,\dots,N_p\subset M$ such that 
        \[
        \mathfrak t_+^{{\rm ss}}=\bigcap_{i=1}^p \{\tau\in \mathfrak t_+\mid \mu(N_{i,\tau})\leq \mu(M_\tau)\}.
        \] 
        \item There is a canonical decomposition of $\mathfrak t_+$ into finitely many chambers $\sigma_1,\dots,\sigma_K$ such that 
        \begin{itemize}
            \item For each $\sigma_i$, there is a filtration $F^i_\bullet: 0\subsetneq M^i_1\subsetneq M^i_2\subsetneq \dots\subsetneq M^i_{l_i}=M$ by saturated $\mathbb X^*(T_\mathbb C)$-graded submodules such that $(F^i_{\bullet})_\tau$ gives the Harder-Narasimhan filtration of $M_\tau$ for every $\tau\in \sigma_i$.
            \item If $i\neq j$, the corresponding Harder-Narasimhan filtrations $F^i_\bullet$ and $F^j_\bullet$ are different.
        \end{itemize}
    \end{enumerate}
\end{algtheorem}
\begin{defn}
    We call the chamber decomposition $\Sigma=\{\sigma_1,\dots,\sigma_K\}$ satisfying the condition in (iii) {\bf HN chamber decomposition}. In this terminology, chambers are not required to be open or connected.
\end{defn}
\begin{proof}
    For any saturated $\mathbb X^*(T_\mathbb C)$-graded submodule $N\subset M$ of rank $s$, we have $f\in F_s^{\max}$ such that $\mu(N)\leq f$. Indeed, the argument in the proof of Theorem \ref{finiteness of maximal slope functions} shows that for any $\mu(N)\in F_{s,i}$, there exists $g\in F^{\max}_{s,i}$ such that $\mu(N)\leq g$. Since $\bigcup_i F^{\max}_{s,i}$ is a finite set, there exists $f\in F^{\max}_s$ with $g\leq f$. 
    Thus, by Theorem \ref{finiteness of maximal slope functions}, we have 
    \[
    \mathfrak t_+^{\st}(M)=\bigcap_{f\in \bigcup_s F^{\max}_s}\{f<\mu(M)\},
    \]
    where $\bigcup_s F^{\max}_s$ is a finite set. This proves (i).
    
    Let $\mathfrak F^{\max}$ be the collection of all saturated graded submodules of $M$ which occur as maximal destabilizing submodules for some $\tau\in \mathfrak t_+$. By the uniqueness of maximal destabilizing submodules, we have $\mu(N)\in F^{\max}_{\rank(N)}$ for any $N\in \mathfrak F^{\max}$, and the map
    \[
    \mu: \mathfrak F^{\max}\to \bigcup_{s=1}^{\rank(M)-1}F_s^{\max}
    \] 
    is injective, hence $\mathfrak F^{\max}$ is finite by Theorem \ref{finiteness of maximal slope functions}. 

    Let $\mathfrak F^{\max}=\{N_1,\dots,N_p\}$. Then we have 
    \[
    \mathfrak t_+^{{\rm ss}}(M)=\bigcap_{i=1}^p\{\mu(N_i)\leq \mu(M)\}.
    \]
    In particular, we have (ii) of the statement. Let $C_i$ be the cone defined as follows: 
    \[
    C_i=\{\mu(N_i)> \mu(M)\}\cap \bigcap_{N_j\supsetneq N_i}\{\mu(N_j)<\mu(N_i)\}\cap \bigcap_{j\neq i}\{\mu(N_j)\leq \mu(N_i)\}.
    \]
    Then we have 
    \[
    \mathfrak t_+=\mathfrak t^{{\rm ss}}(M)\sqcup \bigsqcup_{i=1}^p C_i.
    \]
    For each $C_i$ and any $\tau\in C_i$, $N_{i,\tau}$ is the maximal destabilizing submodule of $M_\tau$.

    For each $i$, let $\mathfrak F^2_i$ be the collection of saturated $\mathbb X^*(T_\mathbb C)$-graded submodules of $M/N_i$ which defines maximal destabilizing submodule of $M/N_i$ for some $\tau\in C_i$. By the same reasoning as above, $\mathfrak F^2_i$ is finite and we can decompose $C_i$ into the cones
    \[
    C_i=(\mathfrak t_+^{{\rm ss}}(M/N_i)\cap C_i)\sqcup \bigsqcup_{j=1}^q C_{ij},
    \]
    which defines the second submodule of the Harder-Narasimhan filtration. 
    Iterating this procedure, we have a decomposition as in (iii). Note that this procedure terminates after finitely many steps since the rank of $M$ is finite.
\end{proof}
\begin{Rem}
Let $\{\sigma_i\}$ be the HN chamber decomposition of $\mathfrak t_+$. If $i\neq j$ and
\[
\sigma_i\cap \bar \sigma_j\neq \emptyset,
\]
then the corresponding filtration $F^i$ is a coarsening of $F^j$. Indeed, by Theorem \ref{finiteness of walls} (ii), as the parameter approaches a boundary point of a chamber, each HN factor remains semistable, while some consecutive HN slopes may become equal. Grouping the consecutive factors having the same limiting slope gives the HN filtration on the limiting stratum.
In particular, we conclude that if $\rank M=2$, then the closures of two distinct unstable HN chambers do not intersect inside the unstable locus.
\end{Rem}
\begin{Rem}
We note that related chamber structures for Harder--Narasimhan filtrations were studied by Neumann \cite[Theorem 3.3.4, Proposition 3.3.5]{neumann2009decomposition} in the setting of stability of sheaves on projective manifolds with respect to movable curve classes. In particular, he proved that the corresponding chamber decomposition is locally finite in the interior of the movable cone, and finite when the movable cone is polyhedral. He also noted that the HN chamber decomposition of the tangent sheaf on the blow-up of $\mathbb P^2$ at the nine base points of a pencil generated by two general cubic curves has infinitely many HN chambers.
In contrast, in our setting the parameter is the Reeb vector, and the walls need not be linear when the affine cone is not $\mathbb Q$-factorial.
\end{Rem}
\section{Examples}
\begin{ex}
    Let $R=\mathbb C[z_1,z_2,z_3]$, let an $T_\mathbb C=(\mathbb C^*)^3$ act diagonally on $\mathbb A^3=\Spec(R)$, i.e.
    \[
    (\lambda_1,\lambda_2,\lambda_3)\cdot (z_1,z_2,z_3)=(\lambda_1z_1,\lambda_2z_2,\lambda_3z_3).
    \]
    Let $\chi_i=\deg(z_i)\in \mathbb X^*(T_\mathbb C)$. Consider the $\mathbb X^*(T_\mathbb C)$-graded module $M$ defined by the following exact sequence
    \[
    0\to R(-(\chi_1+\chi_2+\chi_3))\overset{(z_1,z_2,z_3)}\to R(-(\chi_2+\chi_3))\oplus R(-(\chi_1+\chi_3))\oplus R(-(\chi_1+\chi_2))\to M\to 0
    \]
    $M$ is reflexive since $M$ is locally free outside the origin and has projective dimension $1$ at origin. let $\epsilon_i\in R^3$ denote the $i$-th standard basis vector. Then the image of $\epsilon_i$ in $M$, denoted as $e_i$, generates a saturated submodule of $M$. Indeed, we have an isomorphism 
    \[
    M/Re_i\simeq \sum_{j\neq i}Rz_j\subset R
    \]
    as an ungraded $R$ module. Take $\xi\in \mathfrak t_+:=\{\tau\in \mathbb X_*(T_\mathbb C)_{\mathbb R}\mid \langle\chi_i,\tau\rangle>0, \forall i\}$ to be irrational, i.e. $\xi_i$ are $\mathbb Q$-linearly independent if we write 
    \[
    \xi=\sum_{i=1}^3\xi_i\rho_i,
    \]
    where $(\rho_i)\subset \mathbb X_*(T_\mathbb C)$ is dual basis of $(\chi_i)$. Then the Reeb flow $\phi_t=\exp(\sqrt{-1}\xi t)$ gives a Sasakian structure on $S^5\subset \mathbb C^3$ and restricting $M^\sim$ to $S^5$ gives a structure of $T$-equivariant transverse holomorphic vector bundle $E$ on $(S^5,g,\xi)$. Geometrically, $E$ is pull-back of $T_{\mathbb P^2}(-3)$ by the $T$-equivariant map $\pi: S^5\to \mathbb P^2$, and $\bar\partial_{E,b}$ is induced by the complex structure of $T_{\mathbb P^2}(-3)$.
    
    Take $\tau=\sum_i \tau_i\rho_i \in \mathfrak t_+$. Let $L\subset M_\tau$ be a saturated rank $1$ graded submodule. Then $L$ is a free module generated by some homogeneous element $l\in M_{\tau,\lambda}$ since $L$ is reflexive and $R$ is factorial. $l$ can be written as a $R_\tau$-linear combination of $e_i$'s, and hence we have 
    \[
    \deg_\tau(l)-\deg_\tau(e_i)\geq 0
    \]
    for some $i$. Since
    \[
    \mu_{\gr}(L)=-\lambda\leq -\langle\deg(e_i),\tau\rangle =\mu_{\gr}((Re_i)_\tau),
    \]
    for some $i$, the stability of $M_\tau$ can be tested by considering the slope of $(Re_i)_\tau$'s. Now we have
    \[
    \begin{cases}
        \mu_{\gr}(M_\tau)=-\frac 1 2(\tau_1+\tau_2+\tau_3),\\
        \mu_{\gr}((Re_1)_\tau)=-(\tau_2+\tau_3),\\
        \mu_{\gr}((Re_2)_\tau)=-(\tau_1+\tau_3),\\
        \mu_{\gr}((Re_3)_\tau)=-(\tau_1+\tau_2).
    \end{cases}
    \]
    Thus $M_\tau$ is semistable if and only if $\tau$ satisfies the following inequalities:
    \[
    \begin{cases}
        \tau_1\leq\tau_2+\tau_3,\\
        \tau_2\leq\tau_1+\tau_3,\\
        \tau_3\leq\tau_1+\tau_2,
    \end{cases}
    \]
    and stable if and only if all of these inequalities are strict. As a consequence, transverse holomorphic vector bundle $E_\tau$ on Sasakian $(S^5,g,\tau)$ is (semi)stable if and only if the above inequalities are satisfied. Moreover, this argument shows that we can decompose $\mathfrak t_+$ by the Harder--Narasimhan filtration of $E_\tau$ as the following figure \ref{HNF-chamber decomposition of Reeb cone by TP2}: The three walls divide the normalized Reeb cone into four open chambers. Equivalently, the Harder–Narasimhan filtration gives a stratification consisting of the semistable central triangle, including its boundary, and the three unstable corner triangles.
    \begin{figure}[htbp]
\centering
\begin{tikzpicture}[scale=1.6, font=\small]

  \coordinate (V1) at (0,{2*sqrt(3)});
  \coordinate (V2) at (-2,0);
  \coordinate (V3) at (2,0);

  \coordinate (M12) at (-1,{sqrt(3)});
  \coordinate (M13) at (1,{sqrt(3)});
  \coordinate (M23) at (0,0);

  \fill[gray!15]
    (V1) -- (M12) -- (M13) -- cycle;

  \fill[gray!15]
    (V2) -- (M12) -- (M23) -- cycle;

  \fill[gray!15]
    (V3) -- (M13) -- (M23) -- cycle;

  \fill[gray!35]
    (M12) -- (M13) -- (M23) -- cycle;

  \draw[very thick]
    (V1) -- (V2) -- (V3) -- cycle;

  \draw[thick]
    (M12) -- (M13)
    (M12) -- (M23)
    (M13) -- (M23);

  \node[above] at (V1)
    {$\tau_1=1$};

  \node[below left] at (V2)
    {$\tau_2=1$};

  \node[below right] at (V3)
    {$\tau_3=1$};

  \node[rotate=60, right] at (-1.3,1.5)
    {$\tau_3=0$};

  \node[rotate=-60, left] at (1.3,1.5)
    {$\tau_2=0$};

  \node[below] at (0,0)
    {$\tau_1=0$};

  \node[align=center] at (0.05,2.1)
    {$\tau_1>\frac12$\\[2pt]
     $\operatorname{HN}_1(M_\tau)=Re_1$};

  \node[align=center] at (-1.0,0.4)
    {$\tau_2>\frac12$\\[2pt]
     $\operatorname{HN}_1(M_\tau)=Re_2$};

  \node[align=center] at (1.0,0.4)
    {$\tau_3>\frac12$\\[2pt]
     $\operatorname{HN}_1(M_\tau)=Re_3$};

  \node[align=center] at (0,1.3)
    {$\mathfrak t^{\st}_+(M)$\\[2pt]
     $\tau_i<\frac12\quad(\forall i)$};

  \node at (0,{2.25*sqrt(3)})
    {$\mathfrak t_+\cap
      \{\tau_1+\tau_2+\tau_3=1\}$};

\end{tikzpicture}
\caption{The Harder--Narasimhan chamber decomposition of the normalized
Reeb cone. The central triangle is the semistable chamber, while the three
corner triangles correspond to the destabilizing submodules
$Re_1,Re_2,Re_3$, respectively. The interior of the central triangle is the stable chamber.}
\label{HNF-chamber decomposition of Reeb cone by TP2}
\end{figure}
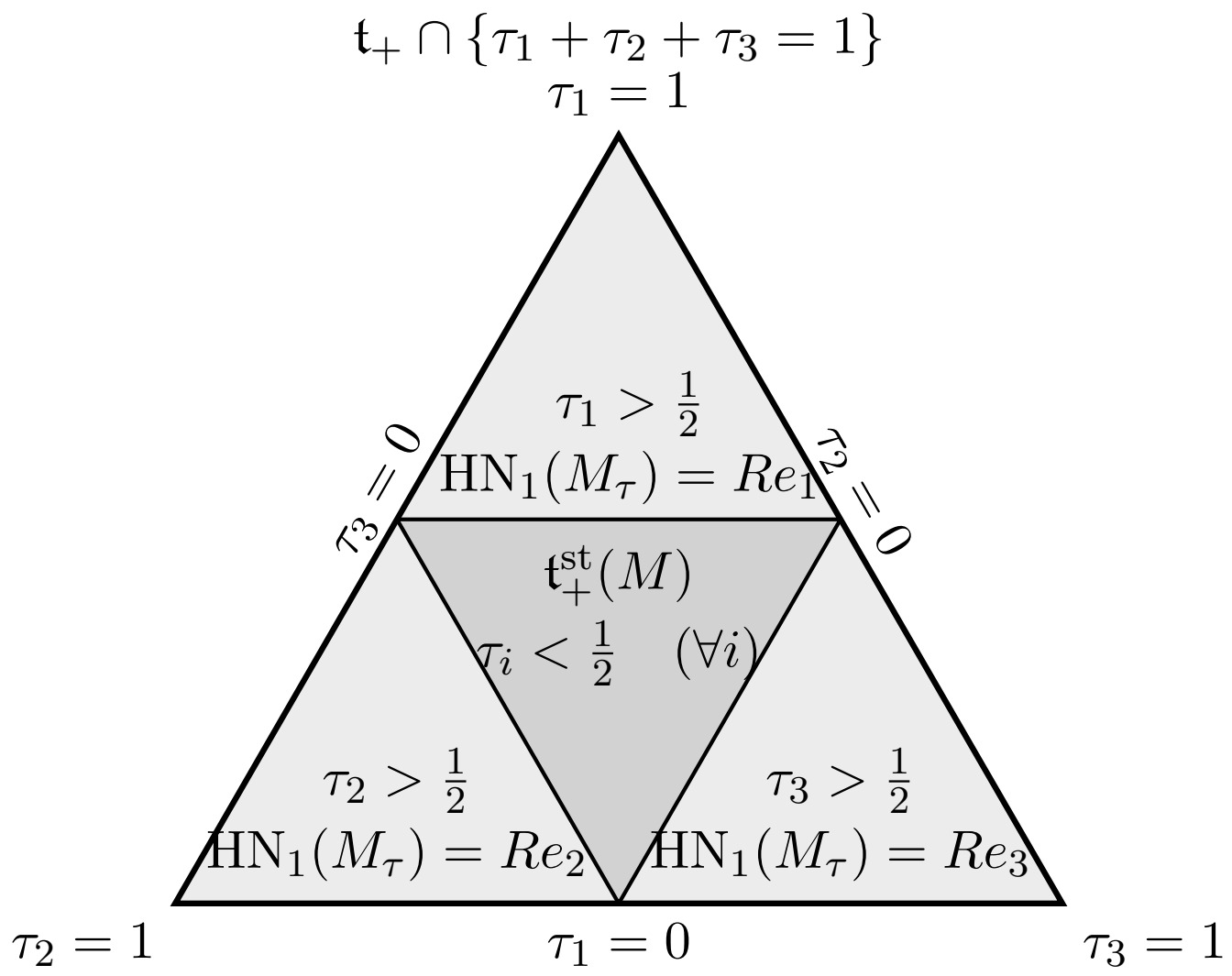

This example yields many simple examples of the following phenomena:
\begin{enumerate}
    \item $\mathbb Q$-Fano variety $X$ with $\rho(X)=1$ whose tangent sheaf is unstable with respect to $-K_X$, and
    \item $\mathbb Q$-Fano variety $X$ with $\rho(X)=1$ whose tangent sheaf is stable, while $X$ is not $K$-semistable.
\end{enumerate}
For (i), choose $(a,b,c)\in \mathbb Z_{>0}$ with $a\leq b\leq c$, $a+b<c$ and $\gcd(a,b)=\gcd(a,c)=\gcd(b,c)=1$. Then $\tau=(a,b,c)$ defines a quasi-regular Reeb vector on $S^5$ and the quotient of $S^5$ by the induced ${\rm U}(1)$-action 
\[
[S^5/{\rm U}(1)_\tau]
\]
has stacky locus with codimension $\geq 2$ and its coarse moduli space is $\mathbb P(a,b,c)$, which is a $\mathbb Q$-Fano surface with $\rho=1$. Since $\mathbb P(a,b,c)$ is normal and the coarse moduli morphism is an isomorphism in codimension $1$, we have
\[
T_{\mathbb P(a,b,c)}\simeq M(\chi_1+\chi_2+\chi_3)_\tau^\sim
\]
and thus it is unstable with respect to $-K_{\mathbb P(a,b,c)}=\mathcal O_{\mathbb P(a,b,c)}(a+b+c)$. In this case, maximal destabilizing subsheaf of $T_{\mathbb P(a,b,c)}$ is $\mathcal O(c)$, corresponding to $R(\chi_1+\chi_2+\chi_3)\cdot e_3\subset M(\chi_1+\chi_2+\chi_3)$, and this defines a foliation on $\mathbb P(a,b,c)$ whose leaves through general points are given by the fibres of the projection 
\[
\mathbb P(a,b,c)\dashrightarrow \mathbb P(a,b),
\]
whose indeterminacy is resolved by the weighted blow-up of $\mathbb P(a,b,c)$ at the point $[0:0:1]$ with weights $(a,b)$. Note that in this case, $[0:0:1]\in \mathbb P(a,b,c)$  is a singular point of type $\frac 1 c(a,b)$, and thus the inequality $a+b<c$ forces the weighted projective space $\mathbb P(a,b,c)$ to have non-canonical singularities.

For (ii), choose any $(1,1,1)\neq (a,b,c)\in \mathbb Z_{>0}$ with $a\leq b\leq c$ and $a+b>c$ and $\gcd(a,b)=\gcd(a,c)=\gcd(b,c)=1$. Then $\tau=(a,b,c)$ defines a quasi-regular Reeb vector on $S^5$. By the same argument as above, the tangent sheaf of $\mathbb P(a,b,c)$ is stable with respect to $-K_{\mathbb P(a,b,c)}$. However, since $(a,b,c)\neq (1,1,1)$, $\mathbb P(a,b,c)$ is not $K$-semistable. Indeed, let $v_1,v_2,v_3$ be primitive generators of rays of corresponding Fan $\Delta$ of $\mathbb P(a,b,c)$. Then we have
\[
av_1+bv_2+cv_3=0.
\]
Let $P$ be the anti-canonical polytope and $p_i$ be the vertex of $P$ which is not in the facet corresponding to $v_i$. Then we have 
\[
\langle p_i,v_j\rangle=-1, \ (i\neq j)
\]
and thus 
\[
\langle p_1,v_1\rangle=\frac{S-a}{a}, \ \langle p_2,v_2\rangle=\frac{S-b}{b}, \ \langle p_3,v_3\rangle=\frac{S-c}{c},
\]
where $S:=a+b+c$. Let $b_P$ be the barycenter of $P$. Then we have 
\[
\langle b_P,v_1\rangle=\frac{S}{3a}-1, \ \langle b_P,v_2\rangle=\frac{S}{3b}-1, \ \langle b_P,v_3\rangle=\frac{S}{3c}-1
\]
and thus $b_P=0$ if and only if $a=b=c$. Since weights are well-formed, this would force $(a,b,c)=(1,1,1)$, contrary to our choice. Therefore $b_P\neq 0$, and the Toric Futaki character is non-zero. In particular, one of the product test configurations associated with the one-parameter subgroup corresponding to the $v_i$ has negative Donaldson-Futaki invariant. Hence $\mathbb P(a,b,c)$ is not $K$-semistable.

We note that Gholampour--Jiang--Kool \cite{GJK17} gave an explicit description of toric sheaves on weighted projective planes, viewed as toric Deligne--Mumford stacks, and studied moduli spaces of stable torsion-free sheaves, in particular in rank two.
\end{ex}
\begin{ex}\label{line bundles on T[1,1]}
    Let $R=\mathbb C[x,y,z,w]/(xw-yz)$. Then $\mathcal C:=\Spec(R)$ is a $3$-dimensional toric variety associated with a cone 
    \[
    \sigma=\langle(0,1,0),(0,0,1),(1,0,-1),(1,-1,0)\rangle_{\mathbb R_{\geq 0}}\subset \mathbb X_*(T_\mathbb C)_{\mathbb R}=\mathfrak t,
    \]
    whose dual cone is given as follows:
    \[
    \omega=\sigma^\vee=\langle(1,0,0),(1,1,0),(1,0,1),(1,1,1)\rangle_{\mathbb R_{\geq 0}}\subset \mathbb X^*(T_\mathbb C)_{\mathbb R}.
    \]
    Hence, the degrees of $x,y,z,w$ are given as follows:
    \[
    \deg(x)=(1,0,0),\ \deg(y)=(1,1,0),\ \deg(z)=(1,0,1),\ \deg(w)=(1,1,1).
    \]
    Let $\tau=(\alpha,\beta,\gamma)\in \mathfrak t_+:={\rm int}(\sigma)$. Then, we can calculate the index character as follows:
    \begin{align*}
        F_R(\tau,t)=&(1-e^{-(2\alpha+\beta+\gamma)t})F_{\mathbb C[x,y,z,w]}(\alpha,\alpha+\beta,\alpha+\gamma,\alpha+\beta+\gamma; t)\\
        =&\frac{1-e^{-(2\alpha+\beta+\gamma)t}}{(1-e^{-\alpha t})(1-e^{-(\alpha+\beta)t})(1-e^{-(\alpha+\gamma)t})(1-e^{-(\alpha+\beta+\gamma)t})}\\
        =&\frac{1}{t^3}\cdot \frac{2\alpha+\beta+\gamma}{\alpha(\alpha+\beta)(\alpha+\gamma)(\alpha+\beta+\gamma)}+O(t^{-2})
    \end{align*}
    Thus we get 
    \[
    \vol(R_\tau)=\frac{2\alpha+\beta+\gamma}{\alpha(\alpha+\beta)(\alpha+\gamma)(\alpha+\beta+\gamma)}.
    \]
    Consider the following homogeneous ideals: 
    \[
    I_1=(x,z),\ I_2=(x,y),\ I_3=(y,w),\ I_4=(z,w)
    \]
    The ideals $I_j$ are non-$\mathbb Q$-Cartier graded reflexive modules of rank $1$. Their slopes are calculated using the exact sequence:
    \[
    0\to I_j\to R\to R/I_j\to 0.
    \]
    For example, if $j=1$, $R/I_1=\mathbb C[y,w]$, and 
    \begin{align*}
        F_{\mathbb C[y,w]}(\alpha+\beta,\alpha+\beta+\gamma;t)=&\frac{1}{(1-e^{-(\alpha+\beta)t})(1-e^{-(\alpha+\beta+\gamma)t})}\\
        =&\frac{1}{t^2}\cdot \frac{1}{(\alpha+\beta)(\alpha+\beta+\gamma)}+O(t^{-1}).
    \end{align*}
    Thus we have 
    \[
    \deg_\tau(I_1)=-\frac{1}{(\alpha+\beta)(\alpha+\beta+\gamma)}.
    \]
    Hence, we conclude that  
    \[
    \mu(I_{1,\tau})=-\frac{\alpha(\alpha+\gamma)}{2\alpha+\beta+\gamma}.
    \]
    Similarly, we have
    \[
    \begin{cases}
        \mu(I_{2,\tau})=-\frac{\alpha(\alpha+\beta)}{2\alpha+\beta+\gamma}\\
        \mu(I_{3,\tau})=-\frac{(\alpha+\beta)(\alpha+\beta+\gamma)}{2\alpha+\beta+\gamma}\\
        \mu(I_{4,\tau})=-\frac{(\alpha+\gamma)(\alpha+\beta+\gamma)}{2\alpha+\beta+\gamma}.
    \end{cases}
    \]
    They are all non-linear with respect to $\tau$. If we intersect the cone $\mathfrak t_+$ by the hyperplane $2\alpha+\beta+\gamma=1$, these level sets are hyperbolic curves.
\end{ex}
\begin{ex}\label{stable cone of some rank 2 bundle on T11}
    Let $\mathcal C$ be as in Example \ref{line bundles on T[1,1]}. Let $e_1,e_2,e_3$ be a standard basis of $\mathbb X^*(T_\mathbb C)$. First, let $\xi_0=(2,0,0)\in \sigma$. Then $\xi_0$ defines a regular Reeb vector corresponding to $\mathcal O_{\mathbb P^1\times \mathbb P^1}(1,1)$. There is a circle bundle 
    \[
    \pi: T^{1,1}:=(SU(2)\times SU(2))/U(1)\simeq S^2\times S^3\to \mathbb P^1\times \mathbb P^1,
    \]
    which can be identified with the unit circle bundle of $\mathcal O(-1,-1)$, i.e. the quotient of the link of $\mathcal C$ to $\mathbb P^1\times \mathbb P^1$. Thus, by pulling back $T$-linearized vector bundles on $\mathbb P^1\times \mathbb P^1$, we have $T$-linearized transverse holomorphic vector bundles on $T^{1,1}$. Now, we have 
    \begin{align*}
        \Ext^1(\mathcal O(1,-1),\mathcal O(-1,0))&=H^1(\mathbb P^1\times \mathbb P^1,\mathcal O(-2,1))\\
        &\simeq H^1(\mathbb P^1,\mathcal O(-2))\otimes H^0(\mathbb P^1,\mathcal O(1))\\
        &\simeq \mathbb C^2.
    \end{align*}
    Now $\mathcal O(1,-1)$ and $\mathcal O(-1,0)$ correspond to graded modules $I_3\cdot I_2^\vee=(y,w)^2(e_1+e_2)$ and $I_2=(x,y)$ respectively. Thus they have torus actions. Fix a character $\chi\in \mathbb X^*(T_\mathbb C)$ with $\langle\xi_0,\chi\rangle=0$ and endow torus actions on $\mathcal O(1,-1)$ and $\mathcal O(-1,0)$, such that whose corresponding modules are
    \[
    (y,w)^2(\chi+e_1+e_2),\ (x,y),
    \]
    respectively. We have 
    \[
    \Ext^1((y,w)^2(\chi+e_1+e_2),(x,y))=\mathbb C_{\chi-e_3}\oplus \mathbb C_{\chi-e_2-e_3}\oplus\mathbb C_{\chi-e_1-e_2-e_3}\oplus \mathbb C_{\chi-e_1-e_2-2e_3}.
    \]
    This computation was performed in Macaulay2 \cite{M2} using the $\mathbb Z^3$-grading
    \[
    \deg(x)=e_1,\ \deg(y)=e_1+e_2,\ \deg(z)=e_1+e_3,\ \deg(w)=e_1+e_2+e_3.
    \]
    The first two direct summands correspond to the basis of $\Ext^1(\mathcal O(1,-1),\mathcal O(-1,0))$ and the last two direct summands correspond to the basis of $\Ext^1(\mathcal O(2,0),\mathcal O(-1,0))$. Now, since $\Hom(\mathcal O(-1,0),\mathcal O(1,-1))=H^0(\mathbb P^1\times \mathbb P^1,\mathcal O(2,-1))=0$, no non-trivial extension has middle term isomorphic to $\mathcal O(-1,0)\oplus \mathcal O(1,-1)$.

    Fix $\chi=e_3$. Then we can choose a $T$-linearized non-trivial extension
    \[
    0\to \mathcal O(-1,0)\to \mathcal E\to \mathcal O(1,-1)\to 0.
    \]
    of vector bundles. $\mathcal E$ is stable sheaf on $\mathbb P^1\times \mathbb P^1$ with respect to $\mathcal O(1,1)$. In fact, let $\mathcal L\subset \mathcal E$ be a saturated subsheaf of $E$. Since $\mathcal L$ is reflexive sheaf on smooth surface, $\mathcal L$ is isomorphic to some line bundle $\mathcal O(a,b)$. We have 
    \[
    \mu(\mathcal E)=-\frac 12.
    \]
    Thus, if $\mathcal L$ is destabilising, we have 
    \[
    a+b\geq 0,
    \]
    Thus $\mathcal L\not\subset \mathcal O(-1,0)$. 
    
    Let $\mathcal L\to \mathcal E\to \mathcal O(1,-1)$ be a composition of morphism $\mathcal L\subset \mathcal E$ and projection $\mathcal E\to \mathcal O(1,-1)$, which is non-zero since $\mathcal L\not\subset \mathcal O(-1,0)$. Since $\mathcal L$ is rank $1$ torsion-free sheaf, this morphism is injective. This shows that 
    \[
    a\leq 1,\ b\leq -1.
    \]
    Thus we have $a=1$ and $b=-1$. This means that the morphism $\mathcal E\to \mathcal O(1,-1)$ splits, contradicting our choice of $\mathcal E$. Thus $\mathcal E$ is stable. Thus corresponding $\mathbb X^*(T_\mathbb C)$-graded module $M=M(\mathcal E)$ is also stable with respect to $\xi_0=(2,0,0)$. 

    Let $N\subset M$ be a saturated rank $1$ $\mathbb X^*(T_\mathbb C)$-graded submodule. By the exact sequence
    \[
    0\to I_2\to M\to I_3^2(e_1+e_2+e_3)\to 0,
    \]
    we have either $N\subset I_2$ or $N$ injects to $I_3^2(e_1+e_2+e_3)$. In fact, if $N\not\subset I_2$, we have an injection
    \[
    I_2/(N\cap I_2)\to M/N.
    \]
    Since $M/N$ is torsion-free, $I_2/(N\cap I_2)$ is either $0$ or rank $1$ graded submodule of $M/N$. Since both $N$ and $I_2$ are torsion-free of rank $1$ and $N\not\subset I_2$, we have $N\cap I_2=0$. Thus $N$ injects to $I_3^2(e_1+e_2+e_3)$. If $N\subset I_2$, we have 
    \[
    \mu_\tau(N)\leq \mu_\tau(I_2)
    \]
    for any $\tau\in \mathfrak t_+$. 

    Let $N\not\subset I_2$, then $N_{\xi_0}$ induces a subsheaf $\mathcal O(a,b)\subset \mathcal O(1,-1)$. Thus we have 
    \[
    a\leq 1, \ b\leq -1.
    \]
    Since $N\subset M$, the inclusion $s: \mathcal O(a,b)\subset \mathcal O(1,-1)$ comes from $\tilde s: \mathcal O(a,b)\to \mathcal E$ and thus $s$ maps to $0$ in 
    \[
    \Ext^1(\mathcal O(a,b),\mathcal O(-1,0))\simeq H^1(\mathbb P^1,\mathcal O(-1-a))\otimes H^0(\mathbb P^1,\mathcal O(-b)).
    \]
    Let $e\in \Ext^1(\mathcal O(1,-1),\mathcal O(-1,0))$ be an extension class corresponding to $\mathcal E$. Then the image of $s$ in $\Ext^1(\mathcal O(a,b),\mathcal O(-1,0))$ is exactly $e\cdot s$. Since $e\neq 0$, existence of $\tilde s$ is equivalent to $e\cdot s=0$. If $a\leq 0$, we have 
    \[
    H^1(\mathbb P^1,\mathcal O(-1-a))\otimes H^0(\mathbb P^1,\mathcal O(-b))=0.
    \]
    If $a=1$, $s: H^1(\mathbb P^1,\mathcal O(-2))\otimes H^0(\mathbb P^1,\mathcal O(1))\to H^1(\mathbb P^1,\mathcal O(-2))\otimes H^0(\mathbb P^1,\mathcal O(-b))$ is injective. Hence $e\cdot s\neq 0$. Thus, the existence of $\tilde s$ is equivalent to $a\leq 0$, $b\leq -1$.

    For $(a,b)=(0,-1)$, we have two $\mathbb X^*(T_\mathbb C)$-graded morphisms
    \[
    y: I_3(e_3)\to I_3^2(e_1+e_2+e_3),
    \]
    and
    \[
    w: I_3\to I_3^2(e_1+e_2+e_3).
    \]
    These morphisms satisfy the following properties: for any saturated graded submodule $N\subset M$ with $N\not\subset I_2$, the induced injection $N\subset I_3^2(e_1+e_2+e_3)$ factors through one of the above injections. Indeed, for such $N$, we have graded injection $N\to I_3^2(e_1+e_2+e_3)$ and this induces a $T$-invariant element $s\in H^0(\mathbb P^1\times \mathbb P^1,\mathcal O(1-a,-1-b))$, with $(a,b)$ corresponding to $N_{\xi_0}$. Thus $s$ can be written as 
    \[
    s=\lambda X_1^iX_2^{1-a-i}Y_1^jY_2^{-1-b-j},
    \]
    where $X_1$, $X_2$ are the homogeneous coordinates of the first component, $Y_1$, $Y_2$ are the homogeneous coordinates of the second component, $\lambda\in \mathbb C^*$, $0\leq i\leq 1-a$ and $0\leq j\leq -1-b$. 
    Since $1-a\geq 1$, $s: \mathcal O(a,b)\to \mathcal O(1,-1)$ factors through 
    \[
    \mathcal O(a,b)\to \mathcal O(0,-1)\overset{X_1}\to \mathcal O(1,-1),
    \]
    or
    \[
    \mathcal O(a,b)\to \mathcal O(0,-1)\overset {X_2}{\to}\mathcal O(1,-1).
    \]
    These two morphisms correspond to $y: I_3(e_3)\to I_3^2(e_1+e_2+e_3)$ and $w: I_3\to I_3^2(e_1+e_2+e_3)$ respectively.
    Thus the injection $N\to I_3^2(e_1+e_2+e_3)$ factors through one of these two morphisms and thus the claim holds. As we have seen above, we have lifts of $y: I_3(e_3)\to I_3^2(e_1+e_2+e_3)$ and $w: I_3\to I_3^2(e_1+e_2+e_3)$ to $\tilde y: I_3(e_3)\to M$ and $\tilde w: I_3\to M$. By maximality of these submodules, images of $\tilde y$ and $\tilde w$ are saturated. 

    Thus, for any saturated rank $1$ graded submodule $N$ of $M$, we have 
    \begin{itemize}
        \item Either $N\subset I_2$, or 
        \item $N\hookrightarrow I_3$, or
        \item $N\hookrightarrow I_3(e_3)$.
    \end{itemize}
    This shows that 
    \[
    \mu_\tau(N)\leq \max\{\mu_\tau(I_2),\mu_\tau(I_3),\mu_\tau(I_3(e_3))\}.
    \]
    Thus we have the following proposition: 
    \begin{prop}
        Let $\tau\in \mathfrak t_+$. If $M_\tau$ is unstable, the maximal destabilizing submodule of $M_\tau$ must be isomorphic to one of the following $\mathbb X^*(T_\mathbb C)$-graded modules:
        \[
        (I_2)_\tau, \ (I_3)_\tau,\ I_3(e_3)_\tau
        \]\qed
    \end{prop}
    Set $\tau=(\alpha,\beta,\gamma)$. Then we have 
    \[
    \begin{cases}
        \mu_\tau(I_2)=-\frac{\alpha(\alpha+\beta)}{2\alpha+\beta+\gamma},\\
        \mu_\tau(I_3)=-\frac{(\alpha+\beta)(\alpha+\beta+\gamma)}{2\alpha+\beta+\gamma},\\
        \mu_\tau(I_3(e_3))=\frac{(\alpha+\gamma)\gamma-(\alpha+\beta)^2}{2\alpha+\beta+\gamma}
    \end{cases}
    \]
    and 
    \[
    \mu_\tau(M)=\frac{1}{2}\mu(I_3(e_3))=\frac{(\alpha+\gamma)\gamma-(\alpha+\beta)^2}{2(2\alpha+\beta+\gamma)}.
    \]
    Thus the stable locus is defined by the following inequalities: 
    \[
    \begin{cases}
        \alpha^2-\beta^2+\gamma^2+\alpha\gamma>0\\
        \alpha^2+\beta^2+\gamma^2+2\alpha\beta+2\beta\gamma+3\alpha\gamma>0\\
        \alpha^2+\beta^2-\gamma^2+2\alpha\beta-\alpha\gamma>0.
    \end{cases}
    \]
    If we set $u=\alpha+\beta$, $v=\alpha+\gamma$, cutting by the equation $u+v=1$, the stable locus and chamber decomposition of $\mathfrak t_+$ with respect to maximal destabilizing subsheaves can be drawn as the figure \ref{fig:HN-chamber-M}.

The same argument shows that if we choose $\chi=e_2+e_3$, we have 
\[
\mu_\tau(N)\leq \max\{\mu_\tau(I_2),\mu_\tau(I_1),\mu_\tau(I_1(e_3))\}
\]
and 
\[
\begin{cases}
\mu_\tau(I_2)=-\frac{\alpha(\alpha+\beta)}{2\alpha+\beta+\gamma}\\
\mu_\tau(I_1)=-\frac{\alpha(\alpha+\gamma)}{2\alpha+\beta+\gamma}\\
\mu_\tau(I_1(e_3))=\frac{(\alpha+\beta+\gamma)\gamma-\alpha^2}{2\alpha+\beta+\gamma}\\
\mu_\tau(M')=\frac{1}{2}\mu(I_1(e_3))=\frac{(\alpha+\beta+\gamma)\gamma-\alpha^2}{2(2\alpha+\beta+\gamma)},
\end{cases}
\]
where $M'$ is the corresponding $\mathbb X^*(T_\mathbb C)$-graded module. 

Thus, if we intersect the Reeb cone with the hyperplane $2\alpha+\beta+\gamma=1$, we obtain the chamber decomposition shown in Figure \ref{fig:HN-chamber-Mprime}.

Note that the two chamber decompositions are related by the reflection accross the diagonal $\alpha=u$. This symmetry can be explained as follows. We have a map
\[
R=\mathbb C[x,y,z,w]/(xw-yz)\simeq  \mathbb C[ac,ad,bc,bd]\subset \mathbb C[a,b,c,d]
\]
defined by 
\[
\begin{pmatrix}
x&y\\
z&w\\
\end{pmatrix}
=\begin{pmatrix}
a\\
b\\
\end{pmatrix}\cdot 
\begin{pmatrix}
c&d\\
\end{pmatrix}.
\]
This isomorphism corresponds to the Segre embedding $\mathbb P^1_{[a:b]}\times \mathbb P^1_{[c:d]}\to \mathbb P^3$. We endow $\mathbb C[a,b,c,d]$ with the following $\mathbb Z^3$-grading, which is compatible with that of $R$;
\[
\wt(a)=e_1-e_3,\ \wt(b)=e_1,\ \wt(c)=e_3,\ \wt(d)=e_2+e_3.
\]
Now, $\mathbb P^1$ has a natural $\SL(2,\mathbb C)$-action, and thus the maximal compact subgroup $\SU(2)\subset \SL(2,\mathbb C)$ also acts on $\mathbb P^1$. This shows that $\SU(2)\times \SU(2)$ acts on $\mathbb P^1_{[a:b]}\times \mathbb P^1_{[c:d]}$, where the first component acts on $\mathbb P^1_{[a:b]}$ and second on $\mathbb P^1_{[c:d]}$. This action induces actions on $\mathcal O(p,q)$ for any $p,q\in \mathbb Z$. In particular, there is an induced $\SU(2)\times \SU(2)$-action on
\[
\Ext^1(\mathcal O(1,-1),\mathcal O(-1,0))\simeq H^1(\mathbb P^1,\mathcal O(-2))\otimes H^0(\mathbb P^1,\mathcal O(1)).
\]
By Serre duality
\[
H^1(\mathbb P^1,\mathcal O(-2))=H^0(\mathbb P^1,\mathcal O)^*\simeq \mathbb C
\]
is the trivial representation of the first $\SU(2)$, whereas $H^0(\mathbb P^1,\mathcal O(1))=\mathbb C\cdot c+\mathbb C\cdot d$ is the standard representation of the second $\SU(2)$. Hence, as an $\SU(2)\times \SU(2)$-representation, we have 
\[
\Ext^1(\mathcal O(1,-1),\mathcal O(-1,0))\simeq \mathbf{1}\boxtimes \mathbf{2}.
\]
If we restrict this representation to ${\rm U}(1)\subset \{1\}\times \SU(2)$, we have a weight decomposition
\[
H^1(\mathbb P^1,\mathcal O(-2))\otimes H^0(\mathbb P^1,\mathcal O(1))=\mathbb C \cdot c\oplus \mathbb C\cdot d,
\]
where the weight of $c$ is $+1$ and that of $d$ is $-1$. Let 
\[
\kappa=\begin{pmatrix}
0&-1\\
1&0\\
\end{pmatrix}\in \SU(2),
\]
which represents a non-trivial element of the Weyl group $W(\SU(2))\simeq \mathbb Z/2$. In the induced action on coordinate functions, $\kappa$ acts $c\mapsto d$ and $d\mapsto -c$, especially, $\kappa$ flips two weight vectors on the weight decomposition. 

By definition of $\kappa$, we have 
\[
\kappa\cdot x=y, \ \kappa\cdot y=-x, \ \kappa\cdot z=w, \ \kappa\cdot w=-z.
\]
Thus, $\kappa$ defines
\[
\kappa: I_1\mapsto I_3, \ I_2\mapsto I_2,\ I_3\mapsto I_1,\ I_4\mapsto I_4.
\]
With respect to the grading, $\kappa$ maps to 
\[
e_1\mapsto e_1+e_2, \ e_2\mapsto -e_2,\ e_3\mapsto e_3,
\]
especially, we have 
\[
I_2\mapsto I_2, \ I_3\mapsto I_1, \ I_3(e_3)\mapsto I_1(e_3).
\]
On $\mathfrak t_+$, we have 
\[
\kappa: \mathfrak t_+\ni (\alpha,\beta,\gamma)\mapsto (\alpha+\beta,-\beta,\gamma)\in \mathfrak t_+,
\]
If $u=\alpha+\beta$, then the induced action satisfies
\[
(\alpha,u)\mapsto (u,\alpha)
\]
Hence, on the normalized slice, the Weyl reflection is precisely the reflection across the diagonal $\alpha=u$, and it identifies the chamber decompositions in Figure \ref{fig:HN-chamber-both}.

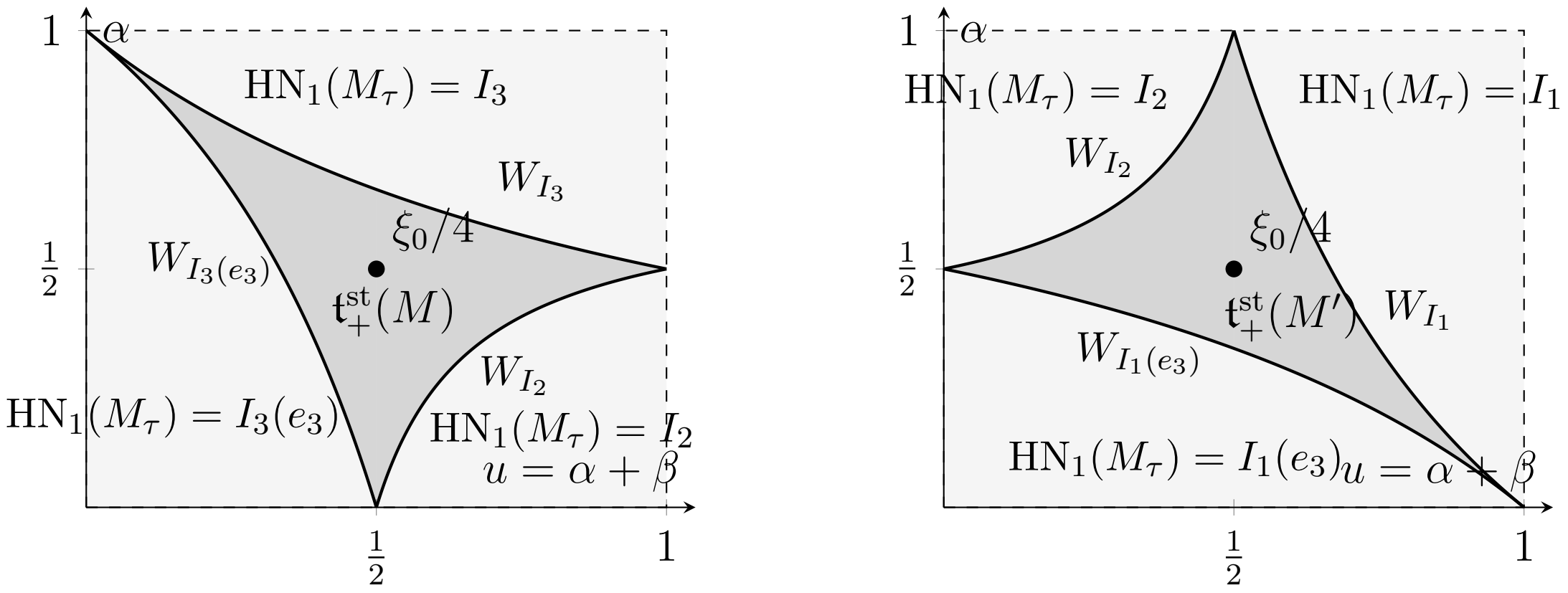
\begin{figure}[htbp]
\centering

\begin{subfigure}[t]{0.48\textwidth}
\centering
\begin{tikzpicture}
\begin{axis}[
    width=7.2cm,
    height=6.2cm,
    axis lines=middle,
    xmin=0, xmax=1.05,
    ymin=0, ymax=1.05,
    xlabel={$u=\alpha+\beta$},
    ylabel={$\alpha$},
    xtick={0,0.5,1},
    xticklabels={$0$,$\frac12$,$1$},
    ytick={0,0.5,1},
    yticklabels={$0$,$\frac12$,$1$},
    samples=200,
    clip=false,
]

\draw[dashed]
    (axis cs:0,0) -- (axis cs:1,0)
    -- (axis cs:1,1) -- (axis cs:0,1) -- cycle;

\addplot[
    domain=0:1,
    thick,
    name path=upper
]
{1/(1+x)};

\addplot[
    domain=0:0.5,
    thick,
    name path=lowerleft
]
{(1-2*x)/(1-x)};

\addplot[
    domain=0.5:1,
    thick,
    name path=lowerright
]
{(2*x-1)/(3*x-1)};

\addplot[
    domain=0:0.5,
    draw=none,
    name path=upperleft
]
{1/(1+x)};

\addplot[
    domain=0.5:1,
    draw=none,
    name path=upperright
]
{1/(1+x)};

\addplot[
    domain=0:1,
    draw=none,
    name path=top
]
{1};

\addplot[
    domain=0:0.5,
    draw=none,
    name path=bottomleft
]
{0};

\addplot[
    domain=0.5:1,
    draw=none,
    name path=bottomright
]
{0};

\addplot[
    fill=gray,
    fill opacity=0.08,
    draw=none
]
fill between[
    of=top and upper
];

\addplot[
    fill=gray,
    fill opacity=0.08,
    draw=none
]
fill between[
    of=lowerleft and bottomleft
];

\addplot[
    fill=gray,
    fill opacity=0.08,
    draw=none
]
fill between[
    of=lowerright and bottomright
];

\addplot[
    fill=gray,
    fill opacity=0.32,
    draw=none
]
fill between[
    of=upperleft and lowerleft
];

\addplot[
    fill=gray,
    fill opacity=0.32,
    draw=none
]
fill between[
    of=upperright and lowerright
];

\addplot[
    only marks,
    mark=*,
    mark size=2pt
]
coordinates {(0.5,0.5)};

\node[anchor=south west]
    at (axis cs:0.5,0.5)
    {$\xi_0/4$};

\node[anchor=south west, font=\small]
    at (axis cs:0.68,0.61)
    {$W_{I_3}$};

\node[anchor=south east, font=\small]
    at (axis cs:0.35,0.43)
    {$W_{I_3(e_3)}$};

\node[anchor=north west, font=\small]
    at (axis cs:0.65,0.35)
    {$W_{I_2}$};

\node[align=center, font=\small]
    at (axis cs:0.50,0.88)
    {$\operatorname{HN}_1(M_\tau)=I_3$};

\node[align=center, font=\small]
    at (axis cs:0.15,0.19)
    {$\operatorname{HN}_1(M_\tau)=I_3(e_3)$};

\node[align=center, font=\small]
    at (axis cs:0.82,0.16)
    {$\operatorname{HN}_1(M_\tau)=I_2$};

\node[align=center]
    at (axis cs:0.53,0.41)
    {$\mathfrak t^{\st}_+(M)$};

\end{axis}
\end{tikzpicture}
\caption{The chamber decomposition for $M$.}
\label{fig:HN-chamber-M}
\end{subfigure}
\hfill
\begin{subfigure}[t]{0.48\textwidth}
\centering
\begin{tikzpicture}
\begin{axis}[
    width=7.2cm,
    height=6.2cm,
    axis lines=middle,
    xmin=0, xmax=1.05,
    ymin=0, ymax=1.05,
    xlabel={$u=\alpha+\beta$},
    ylabel={$\alpha$},
    xtick={0,0.5,1},
    xticklabels={$0$,$\frac12$,$1$},
    ytick={0,0.5,1},
    yticklabels={$0$,$\frac12$,$1$},
    samples=200,
    clip=false,
]

\draw[dashed]
    (axis cs:0,0) -- (axis cs:1,0)
    -- (axis cs:1,1) -- (axis cs:0,1) -- cycle;

\addplot[
    domain=0:0.5,
    thick,
    name path=wallITwo
]
{(1-x)/(2-3*x)};

\addplot[
    domain=0.5:1,
    thick,
    name path=wallIOne
]
{(1-x)/x};

\addplot[
    domain=0:1,
    thick,
    name path=wallIOneTwist
]
{(1-x)/(2-x)};

\addplot[
    domain=0:0.5,
    draw=none,
    name path=topLeft
]
{1};

\addplot[
    domain=0.5:1,
    draw=none,
    name path=topRight
]
{1};

\addplot[
    domain=0:1,
    draw=none,
    name path=bottom
]
{0};

\addplot[
    domain=0:0.5,
    draw=none,
    name path=lowerLeftCopy
]
{(1-x)/(2-x)};

\addplot[
    domain=0.5:1,
    draw=none,
    name path=lowerRightCopy
]
{(1-x)/(2-x)};

\addplot[
    fill=gray,
    fill opacity=0.08,
    draw=none
]
fill between[
    of=topLeft and wallITwo
];

\addplot[
    fill=gray,
    fill opacity=0.08,
    draw=none
]
fill between[
    of=topRight and wallIOne
];

\addplot[
    fill=gray,
    fill opacity=0.08,
    draw=none
]
fill between[
    of=wallIOneTwist and bottom
];

\addplot[
    fill=gray,
    fill opacity=0.32,
    draw=none
]
fill between[
    of=wallITwo and lowerLeftCopy
];

\addplot[
    fill=gray,
    fill opacity=0.32,
    draw=none
]
fill between[
    of=wallIOne and lowerRightCopy
];

\addplot[
    only marks,
    mark=*,
    mark size=2pt
]
coordinates {(0.5,0.5)};

\node[anchor=south west]
    at (axis cs:0.5,0.5)
    {$\xi_0/4$};

\node[anchor=south west, font=\small]
    at (axis cs:0.18,0.66)
    {$W_{I_2}$};

\node[anchor=south west, font=\small]
    at (axis cs:0.73,0.34)
    {$W_{I_1}$};

\node[anchor=north west, font=\small]
    at (axis cs:0.2,0.4)
    {$W_{I_1(e_3)}$};

\node[align=center, font=\small]
    at (axis cs:0.16,0.87)
    {$\operatorname{HN}_1(M_\tau)=I_2$};

\node[align=center, font=\small]
    at (axis cs:0.84,0.87)
    {$\operatorname{HN}_1(M_\tau)=I_1$};

\node[align=center, font=\small]
    at (axis cs:0.40,0.10)
    {$\operatorname{HN}_1(M_\tau)=I_1(e_3)$};

\node[align=center]
    at (axis cs:0.6,0.4)
    {$\mathfrak t^{\st}_+(M')$};

\end{axis}
\end{tikzpicture}
\caption{The chamber decomposition for $M'$.}
\label{fig:HN-chamber-Mprime}
\end{subfigure}

\caption{
Harder--Narasimhan chamber decompositions of the normalized Reeb cone.
Left: the chamber decomposition for $M$, where the maximal destabilizing submodules are $I_3(e_3)$, $I_3$, and $I_2$.
Right: the chamber decomposition for $M'$, where the maximal destabilizing submodules are $I_2$, $I_1$, and $I_1(e_3)$.
In both figures, the darkly shaded region is the stable cone.
}
\label{fig:HN-chamber-both}
\end{figure}

\end{ex}
\section*{Acknowledgements}
This project grew out of discussions with Xuemiao Chen during my stay in Canada. I am grateful to him for many stimulating conversations that led to the development of this work. 
I would also like to thank my advisor, Yuji Odaka, for his guidance and many helpful comments on the manuscript. 
I am grateful to Masafumi Hattori for helpful discussions. In particular, a conversation with him about the possibility of obtaining an openness result in the purely algebraic case motivated me to investigate the questions studied in Section 4.2.
Also, the author would like to thank Keshu Zhou for valuable comments and stimulating discussions related to this work.
\section*{Usage of tools and computational resources}
OpenAI's ChatGPT is used during the preparation of this manuscript for language editing, improving exposition, checking the mathematical arguments by the author, and identifying relevant references. 
In particular, ChatGPT suggested a correction to the argument in the proof of Theorem \ref{Spectral Gap Theorem} and a more direct formulation of the proof of Lemma \ref{restricted volume lemma}. 
These suggestions were subsequently checked and incorporated by the author. Except for them, all other mathematical ideas, the overall strategy and structure of the arguments are the author's own. 
The author takes full responsibility for the content of this paper.

Macaulay2 \cite{M2} was used for the graded $\Ext$ computation in Example \ref{stable cone of some rank 2 bundle on T11}.
\printbibliography
\end{document}